%% file: arxiv_draft_submission.tex
\documentclass[preprint,12pt]{elsarticle}

\usepackage{amsfonts,amsmath,amsthm}
\usepackage{graphicx,color,xcolor} 

\usepackage{graphicx}
\usepackage{subcaption} 

\usepackage{hyperref,float} 

\usepackage{natbib}

\journal{arXiv}

\begin{document}

\begin{frontmatter}

\title{Energy-Dissipative Evolutionary Kolmogorov–Arnold Networks for Complex PDE Systems}


\cortext[cor1]{Corresponding author}
\author[label1,label2,cor1]{Guang Lin}

\author[label2]{Changhong Mou}

\author[label2]{Jiahao Zhang}

\affiliation[label1]{organization={School of Mechanical Engineering, Purdue University},
            addressline={610 Purdue Mall}, 
            city={West Lafayette},
            postcode={47907}, 
            state={IN},
            country={USA}}
\affiliation[label2]{organization={Department of Mathematics, Purdue University},
            addressline={610 Purdue Mall}, 
            city={West Lafayette},
            postcode={47907}, 
            state={IN},
            country={USA}}

\begin{abstract}
 We introduce evolutionary Kolmogorov–Arnold Networks (EvoKAN), a {novel} framework for solving {complex} partial differential equations (PDEs). EvoKAN builds on Kolmogorov–Arnold Networks (KANs), where activation functions are spline-based and trainable on each edge, offering localized flexibility across multiple scales. 
  Rather than retraining the network repeatedly, EvoKAN encodes only the PDE’s initial state during an initial learning phase. The network parameters then evolve numerically, governed by the same PDE, without any additional optimization. 
  By treating these parameters as continuous functions in the relevant coordinates and updating them through time steps, EvoKAN can predict system trajectories over arbitrarily long horizons, a notable challenge for many conventional neural-network-based methods. In addition, EvoKAN integrates the scalar auxiliary variable (SAV) method to guarantee unconditional energy stability and computational efficiency. At individual time step, SAV only needs to solve decoupled linear systems with constant coefficients, the implementation is significantly simplified. We test the proposed framework in several complex PDEs, including one dimensional and two dimensional Allen--Cahn equations and two dimensional Navier-Stokes equations. Numerical results show that EvoKAN solutions closely match analytical references and established numerical benchmarks, effectively capturing both phase-field phenomena (Allen–Cahn) and turbulent flows (Navier–Stokes).

\end{abstract}

\begin{keyword}
Kolmogorov–Arnold Networks (KANs)
\sep 
Scalar Auxiliary Variable (SAV) method
\sep 
Evolutionary Neural Network
\sep 
Partial Differential Equations


\end{keyword}

\end{frontmatter}

\input{main_text}

\end{document}

%% file: main_text.tex
\section{Introduction}
Recently proposed Kolmogorov-Arnold Networks (KANs) \cite{liu2024kan,liu2024kan2}, originating from the Kolmogorov-Arnold Theorem \cite{kolmogorov1957representation,kolmogorov1961representation,braun2009constructive}, have become a promising alternative to traditional multilayer perceptrons (MLPs) in the field of scientific machine learning. Although MLPs have recently received significant attention due to advances in their theory, they are based on trainable weights and biases with predetermined activation functions \cite{apicella2021survey,trentin2001networks}. KANs, on the other hand, use trainable activation functions based on splines, allowing local adjustments and adaptability to different resolutions \cite{liu2024kan,liu2024kan2}. This shift not only enhances the interpretability of KANs over MLPs but also makes them particularly suitable for applications involving continual learning and handling noisy data. 
The applications of KANs have also been explored in the context of solving differential equations. For instance, Liu et al. \cite{liu2024kan} combined KANs with physics-informed neural networks (PINNs) \cite{raissi2019physics} to solve a two-dimensional Poisson equation. Likewise, Abueida et al. introduced DeepOKAN \cite{abueidda2024deepokan}, a KAN operator network based on radial basis functions (RBF), to solve a two-dimensional orthotropic elasticity problem.
Both studies imply that KANs significantly outperformed traditional MLP architectures. However, these works did not explore the potential of KANs in more complex systems.

However, the use of neural networks to approximate solutions to partial differential equations (PDE) has attracted significant interest over the past decade \cite{lu2021learning,li2020fourier,cai2021deepm,guo2021construct,karniadakis2021physics,boulle2023mathematical,kovachki2023neural,zhuang2024two}.
There are two main approaches: The first learns the PDE operator, mapping initial/boundary conditions to solutions using methods such as DeepONet \cite{goswami2022physics,lu2021learning,zhu2023fourier}. Although computationally expensive to train with the available data, these models provide efficient evaluations. The second approach uses neural networks as basis functions to represent a single solution, minimizing residuals and boundary condition mismatches. Among them, physics-informed neural networks (PINNs) \cite{cai2021physics,cuomo2022scientific,toscano2024pinns} have been developed for both forward and inverse problems by enforcing physical laws or information in neural networks. In particular, time-dependent PDEs are solved by minimizing the chosen residuals at pre-selected points across the entire spatio-temporal domain. Although Physics-Informed Neural Networks (PINNs) efficiently compute gradients within PDEs using automatic differentiation, they encounter computational challenges when making long-term predictions of dynamics for time-dependent PDEs due to their lack of temporal causality \cite{beck2021perspective}.

Du et al. \cite{du2021evolutional} introduced Evolution Deep Neural Networks (EDNN), where the network parameters are dynamically updated to predict the PDE's evolution over any time span.  
The parameters in the EDNN are trained solely to capture the initial condition of the system and are subsequently updated recurrently without the need for further training. These updates are performed numerically and are discretized from the governing equations. By advancing the weights in the parameter space, the temporal evolution of the PDE system can be accurately predicted.
The EDNN demonstrates significant potential for approximating complex PDE systems through neural networks.
In this work, we utilize this evolutionary deep operator neural network framework \cite{zhang2024energy} and KANs to design evolutionary KANs \cite{liu2024kan,liu2024kan2} for complex PDE systems. Specifically, we integrate the adaptability and local learning capabilities of KANs with EDNN to create a more robust framework suited for tackling high-dimensional and dynamic problems in physical systems. 
This approach aims to address limitations found in both standard deep neural network models and conventional KAN architectures by replacing the conventional Multi-Layer Perceptrons (MLPs) with KANs that are parametrized as splines. 
By combining the accuracy and interpretability of KANs with the long-term dynamical predictability of EDNN, our proposed model not only achieves improved accuracy, but also enhances generalization to more complex scenarios. 
This combination ensures that the learned models remain consistent with the underlying physical laws in the original governing systems.
In addition, we assess the performance of this framework on various challenging benchmarks, demonstrating its effectiveness in capturing the intricate dynamics of real-world physical phenomena.
In addition, EvoKAN incorporates the scalar auxiliary variable (SAV) methodology \cite{shen2018scalar,shen2019new}, which yields numerical schemes that are unconditionally energy stable and computationally efficient. In particular, these schemes only require solving decoupled linear systems with constant coefficients at each time step, greatly simplifying the implementation.
In brief summary, our primary contribution includes the following:
\begin{enumerate}
\item We propose a novel evolutionary Kolmogorov–Arnold Networks (EvoKAN) for solving complex PDEs. 
\item EvoKAN models Kolmogorov–Arnold network weights as time-dependent functions and updates them through the evolution of the governing PDEs.
\item In addition, we incorporate the scalar auxiliary variable (SAV) method for unconditional energy stability and computational efficiency.
\item Through numerical experiments, we show that our algorithm achieves both high accuracy and rapid convergence.
\end{enumerate}

The remainder of this draft is organized as follows. In Section  \ref{sec:evolution-kan}, we introduce the energy-dissipative evolutionary KAN framework. Section \ref{sec:test} presents numerical tests to evaluate the performance of the proposed framework, utilizing benchmark problems: the one-dimensional and two-dimensional Allen-Cahn equations and the two-dimensional incompressible Navier-Stokes equations. Finally, conclusions and potential future research directions are summarized in Section \ref{sec:conclusion}.

\section{Evolutionary Kolmogorov–Arnold Networks \label{sec:evolution-kan}}

\subsection{Kolmogorov–Arnold Network}
The Kolmogorov–Arnold Network starts from the Kolmogorov–Arnold theorem, which asserts that any multivariate continuous function on a bounded domain can be decomposed into a finite composition of continuous univariate functions and addition operations~\cite{sprecher2002space,koppen2002training,lin1993realization,lai2021kolmogorov,leni2013kolmogorov,fakhoury2022exsplinet,montanelli2020error, he2023optimal}. Specifically, for a smooth function $f: [0,1]^n \to \mathbb{R}$, the following representation holds:
\begin{equation}
    f(\mathbf{x}) = f(x_1,\dots,x_n) = \sum_{q=1}^{2n+1} \Phi_q\left(\sum_{p=1}^n \phi_{q,p}(x_p)\right),\label{eq:KART}
\end{equation}
where $\phi_{q,p}: [0,1] \to \mathbb{R}$ and $\Phi_q: \mathbb{R} \to \mathbb{R}$. This result suggests that only the addition is performed truly on multivariate function, as all others can be expressed using univariate functions and sums. Although this finding seems promising for machine learning, it reduces the problem of learning high-dimensional functions to that of learning a set of one-dimensional functions, which can be non-smooth and thus difficult to learn in practice~\cite{poggio2020theoretical,girosi1989representation}. As a result, the Kolmogorov-Arnold representation theorem, while theoretically valid, has been considered impractical for machine learning applications~\cite{poggio2020theoretical,girosi1989representation}.
\begin{figure}[H]
    \centering
    \includegraphics[width=\linewidth]{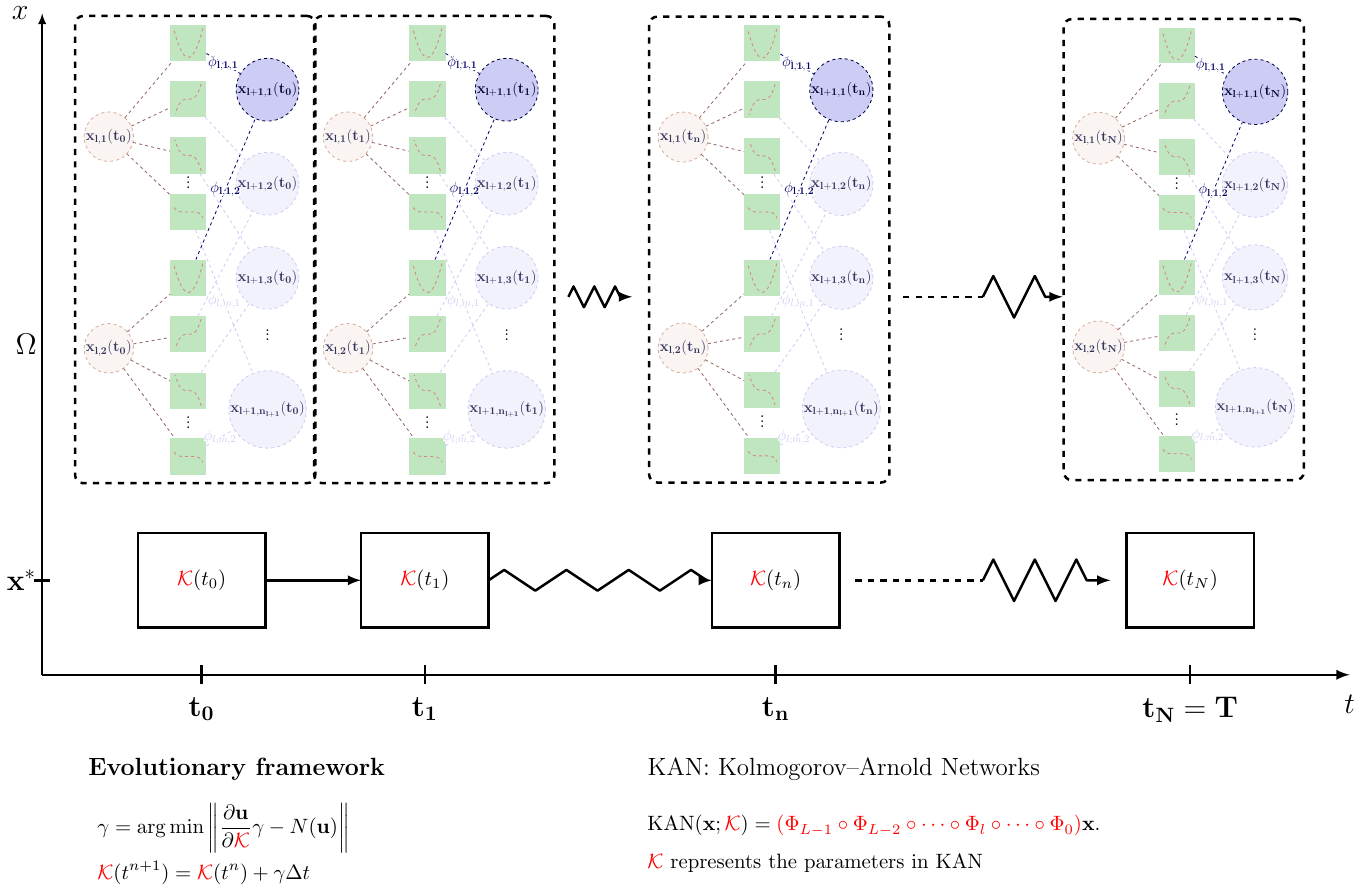}
    \caption{Illustration of Evolutionary Kolmogorov–Arnold Networks (EvoKAN). Kolmogorov–Arnold Networks (KANs) replace linear weights with univariate spline functions learned on each edge. In the evolutionary framework, EvoKAN models Kolmogorov–Arnold Network weights as time-dependent functions and updates them through the evolution of the governing PDEs—without additional training of KAN.
    }
     \label{fig:illustration}
\end{figure}
Given a supervised learning task with input-output pairs $\{\mathbf{x}_i, y_i\}$, the objective is to find a function $f$, i.e., $y_i \approx f(\mathbf{x}_i)$ for each data point. Based on Eq.~(\ref{eq:KART}), this task can be accomplished if certain univariate functions $\phi_{q,p}$ and $\Phi_q$ are identified. This observation motivates the design of a neural network that explicitly parametrizes Eq.~(\ref{eq:KART}). Since the functions are univariate, they can be modeled using B-spline curves with learnable coefficients for the local basis functions. This construction forms the basis for the Kolmogorov-Arnold Network (KAN) architecture, which closely follows the structure of Eq.~(\ref{eq:KART}). The resulting network resembles a two-layer neural network in which activation functions are assigned to edges rather than nodes, and the middle layer has a width of $2n+1$.
Similar to MLPs, where depth is achieved by stacking linear transformations and non-linear activations, a KAN layer in a deeper KAN with $n_{\rm in}$ inputs and $n_{\rm out}$ outputs can be described as a matrix of { single variable} functions:
\begin{equation}
    \mathbf{\Phi} = \{\phi_{q,p}\},\quad p=1,\dots,n_{\rm in},\quad q=1,\dots,n_{\rm out},
\end{equation}
where the functions $\phi_{q,p}$ { contain} trainable parameters. In the Kolmogorov-Arnold theorem, the inner layer corresponds to $n_{\rm in} = n$ and $n_{\rm out} = 2n+1$, while the outer layer corresponds to $n_{\rm in} = 2n+1$ and $n_{\rm out} = 1$. 

The depth of a KAN can be viewed as a series of function compositions, the summation acting as a link between layers. In contrast to MLPs, where activation functions are positioned at the nodes, KANs assign activation functions to the edges. This edge-based configuration simplifies the process of extending KANs in a manner similar to the way that MLPs are deepened. Additionally, KANs naturally incorporate scalar weights into the activation functions, eliminating the need for separate scalar weights. However, in practical implementations, a scalar factor may still be introduced to enhance optimization.

To formalize the structure, a KAN is defined by a sequence of integers:
\begin{equation}
    [n_0, n_1, \dots, n_L],
\end{equation}
where $n_i$ denotes the nodes' number in the $i^{\text{th}}$ layer. The activation value of the $i^{\text{th}}$ neuron in the $l^{\text{th}}$ layer is denoted as $x_{l,i}$. Between layer $l$ and layer $l+1$, $n_l n_{l+1}$ activation functions are existed, with the function connecting neuron $(l,i)$ to neuron $(l+1,j)$ denoted by:
\begin{equation}
    \phi_{l,j,i},\quad l=0,\dots,L-1,\quad i=1,\dots,n_l,\quad j=1,\dots,n_{l+1}.
\end{equation}
The pre-activation of $\phi_{l,j,i}$ is $x_{l,i}$, and the post-activation is denoted as $\tilde{x}_{l,j,i} = \phi_{l,j,i}(x_{l,i})$. The sum of all incoming post-activations consist the activation of neuron $(l+1,j)$:
\begin{equation}
    x_{l+1,j} = \sum_{i=1}^{n_l} \tilde{x}_{l,j,i} = \sum_{i=1}^{n_l} \phi_{l,j,i}(x_{l,i}), \quad j=1,\dots,n_{l+1}.
\end{equation}
This can be written in matrix form as:
\begin{equation}
    \mathbf{x}_{l+1} = 
    \underbrace{\begin{pmatrix}
        \phi_{l,1,1}(\cdot) & \phi_{l,1,2}(\cdot) & \cdots & \phi_{l,1,n_l}(\cdot) \\
        \phi_{l,2,1}(\cdot) & \phi_{l,2,2}(\cdot) & \cdots & \phi_{l,2,n_l}(\cdot) \\
        \vdots & \vdots &  & \vdots \\
        \phi_{l,n_{l+1},1}(\cdot) & \phi_{l,n_{l+1},2}(\cdot) & \cdots & \phi_{l,n_{l+1},n_l}(\cdot) \\
    \end{pmatrix}}_{\mathbf{\Phi}_l}
    \mathbf{x}_l,
\end{equation}
where $\mathbf{\Phi}_l$ is the function matrix {at} the $l^{\text{th}}$ KAN layer. A general KAN network is a composition of $L$ layers, where the output is generated by applying successive KAN layers to an input vector $\mathbf{x}_0 \in \mathbb{R}^{n_0}$:
\begin{equation}\label{eq:kan}
    \text{KAN}(\mathbf{x}) = (\mathbf{\Phi}_{L-1} \circ \mathbf{\Phi}_{L-2} \circ \cdots \circ \mathbf{\Phi}_1 \circ \mathbf{\Phi}_0)\mathbf{x}.
\end{equation}
This equation can also be expressed similarly to Eq.~(\ref{eq:KART}) if the output dimension is $n_L = 1$, defining $f(\mathbf{x}) \equiv \text{KAN}(\mathbf{x})$:
\begin{equation}
    f(\mathbf{x}) = \sum_{i_{L-1}=1}^{n_{L-1}} \phi_{L-1,i_L,i_{L-1}} \cdot
    \left(\sum_{i_{L-2}=1}^{n_{L-2}} \cdots \left(\sum_{i_2=1}^{n_2} \phi_{2,i_3,i_2} \left(\sum_{i_1=1}^{n_1} \phi_{1,i_2,i_1}\left(\sum_{i_0=1}^{n_0} \phi_{0,i_1,i_0}(x_{i_0})\right)\right)\right)\cdots\right).
\end{equation}
Although this expression appears complex, the abstraction of KAN layers provides a clear and intuitive understanding of the structure. The original Kolmogorov-Arnold representation in Eq.~(\ref{eq:KART}) { leads} to a two-layer KAN with shape $[n, 2n+1, 1]$. Since all operations are differentiable, KANs can be trained using backpropagation. In contrast, MLPs are composed of alternating linear transformations $\mathbf{W}$ and nonlinearities $\sigma$:
\begin{equation}
    \text{MLP}(\mathbf{x}) = (\mathbf{W}_{L-1} \circ \sigma \circ \mathbf{W}_{L-2} \circ \sigma \circ \cdots \circ \mathbf{W}_1 \circ \sigma \circ \mathbf{W}_0)\mathbf{x}.
\end{equation}
The primary distinction is that MLPs separate linear transformations and non-linearities, whereas KANs integrate them into function matrices $\mathbf{\Phi}$.

\subsection{Evolutionary Deep Neural Network}
In the following, we present the Evolutionary Deep Neural Network (EDNN) framework \cite{du2021evolutional,zhang2024energy}. In particular, EDNN models neural network weights as time-dependent functions, updating them through the evolution of governing PDEs. Consider a general nonlinear PDE with a initial condition:
\begin{equation}
\label{eqn:general-pde}
 \begin{aligned}
    &\frac{\partial \boldsymbol{u}}{\partial t} + \mathcal{N}_{\boldsymbol{x}}(\boldsymbol{u}) = 0, \\
    &\boldsymbol{u}(\boldsymbol{x}, 0) = \boldsymbol{f}(\boldsymbol{x}),
 \end{aligned}
\end{equation}
where $\boldsymbol{u}(\boldsymbol{x},t) = (u_1, u_2, \dots, u_m)$ is a vector-valued function, $\boldsymbol{x} = (x_1, x_2, \dots, x_d)$ represents the spatial coordinates, and $\mathcal{N}_{\boldsymbol{x}}$ is a nonlinear differential operator.
We can express the approximation of the solution $\bf u$ in $L+1$ layers neural network: 
\begin{equation}
\hat{\boldsymbol{u}}(\boldsymbol{x}, \mathcal{W}(t)) = (\hat{u}_1, \hat{u}_2, \dots, \hat{u}_m)= \mathbf{g}_{L+1}.
\end{equation}
where $\mathbf{g}_{l+1}$ denotes the neurons in the $l+1$-th layer and  $\mathcal{W}(t)$ is the vector that contains all the parameters in the neural network.
A fully connected neural network in a multilayer perceptron (MLP) can be defined as follows:
\begin{equation}\label{NNit}
\mathbf{g}_{l+1}(\mathbf{g}_l) = \sigma(\mathbf{W}_l \mathbf{g}_l + \mathbf{b}_l), \quad l = 0, 1, \dots, L,
\end{equation}
where $\mathbf{g}_l$ denotes the neurons at the $l$-th layer, $\mathbf{W}_l$ and $\mathbf{b}_l$ are the weight matrix and bias vector between layers $l$ and $l+1$, and $\sigma(\cdot)$ is the activation function. 
The network parameters, $\mathbf{W}_l(t)$ and $\mathbf{b}_l(t)$, vary with time, denoted collectively as $\mathcal{W}(t)$.
Moreover, the input to the network is the spatial coordinate of the PDE:
\begin{equation}
 \mathbf{g}_0 = \boldsymbol{x} = (x_1, x_2, \dots, x_d).   
\end{equation}
The time evolution of the solution can be expressed by:
\begin{equation}
\frac{\partial \hat{\boldsymbol{u}}}{\partial t} = \frac{\partial \hat{\boldsymbol{u}}}{\partial \mathcal{W}} \frac{\partial \mathcal{W}}{\partial t}.
\end{equation}
At each time step, the time derivative of the network parameters, $\frac{\partial \mathcal{W}}{\partial t}$, is obtained by solving the optimization problem:
\begin{equation}\label{Wmin}
\frac{\partial \mathcal{W}}{\partial t} = argmin_{\gamma} \mathcal{J}(\gamma), \quad \mathcal{J}(\gamma) = \frac{1}{2} \int_{\Omega} \left\Vert \frac{\partial \hat{\boldsymbol{u}}}{\partial \mathcal{W}} \gamma - \mathcal{N}(\hat{\boldsymbol{u}})\right\Vert_2^2 \mathrm{d}\boldsymbol{x}.
\end{equation}
The optimality condition for this problem is:
\begin{equation}\label{OptCond}
\nabla_{\gamma} \mathcal{J}(\gamma_{\text{opt}}) = \left( \int_{\Omega} \frac{\partial \hat{\boldsymbol{u}}}{\partial \mathcal{W}}^{T} \frac{\partial \hat{\boldsymbol{u}}}{\partial \mathcal{W}} \mathrm{d}\boldsymbol{x} \right)\gamma_{\text{opt}} - \left( \int_{\Omega} \frac{\partial \hat{\boldsymbol{u}}}{\partial \mathcal{W}}^{T} \mathcal{N}(\hat{\boldsymbol{u}}) \mathrm{d}\boldsymbol{x} \right) = 0.
\end{equation}
The optimal $\gamma_{\text{opt}}$ is approximated by:
\begin{equation}\label{OptCondApprox}
\mathbf{J}^{T} \mathbf{J} \hat{\gamma}_{\text{opt}} = \mathbf{J}^{T} \mathbf{N},
\end{equation}
where $\mathbf{J}$ is the network's Jacobian and $\mathbf{N}$ is the PDE residual evaluated at collocation points:
\begin{equation}
\left( \mathbf{J} \right)_{ij} = \frac{\partial \hat{\boldsymbol{u}}^i}{\partial \mathcal{W}_j}, \quad \left( \mathbf{N} \right)_i = \mathcal{N}(\hat{\boldsymbol{u}}^i),
\end{equation}
where $i = 1, 2, \dots, N_u$ indexes the collocation points, and $j = 1, 2, \dots, N_{\mathcal{W}}$ indexes the network parameters. Automatic differentiation is used to compute the elements of $\mathbf{J}$ and $\mathbf{N}$.
We can therefore update the weights, $\mathcal{W}^{n+1}$ and the solution, $\hat{\boldsymbol{u}}^{n+1}(\boldsymbol{x}, \mathcal{W}(t))$ at $n+1$ time step as:
\begin{align}
\mathcal{W}^{n+1} &=\mathcal{W}^n +\hat{\gamma}_{\text{opt}}  \Delta t, \\
\hat{\boldsymbol{u}}^{n+1}(\boldsymbol{x}, \mathcal{W}(t))&=\hat{\boldsymbol{u}}(\boldsymbol{x}, \mathcal{W}^{n+1}(t)).
\end{align}
\subsection{Scalar Auxiliary Variable (SAV) Approach}
\cite{shen2018scalar,shen2019new} introduced the scalar auxiliary variable (SAV) approach, which retains all the advantages of the invariant energy quadratization (IEQ) method while avoiding most of its limitations. The advantages of SAV are two-fold: (1). For single-component gradient flows, at each time step, it leads to linear equations with constant coefficients, making the implementation remarkably straightforward.
(2). For multicomponent gradient flows, it results in decoupled linear equations, one for each component, each with constant coefficients, further simplifying the computational procedure.

In fact, at each time step, the method requires solving a linear system in the form 
\[
(1 - c_0 \Delta t \,\mathcal{G}\mathcal{L})\bar{x} = \bar{b}
\]
twice, where $c_0$ is a positive constant dependent on the time-discretization scheme. Hence, solving the gradient flow using the SAV approach is equivalent to solving the linear parabolic PDE 
\[
\frac{\partial \phi}{\partial t} = \mathcal{G}\mathcal{L} \phi
\]
with an implicit scheme. When $\mathcal{G}$ and $\mathcal{L}$ commute, the linear system is symmetric; even if they do not, the system matrix remains independent of time.
Furthermore, the SAV approach only requires the integrated free energy 
\[
E_1[\phi]=\int_\Omega g(\phi)\,d{\bf x}
\]
to be bounded from below, rather than imposing a uniform lower bound on the free energy density $g(\phi)$, allowing us to handle a broader class of free energies.

We consider equation \eqref{Gflow} with the free energy in the form such that
$\mathcal{E}_1[\phi]$ is bounded from below:
\begin{equation}\label{Gflow}
  \frac{\partial \phi}{\partial t}=\mathcal{G}\mu, 
\end{equation}
{To derive the general form of SAV}, we assume that $\mathcal{E}_1[\phi]\ge C_0>0$. {This assumption is reasonable provided that if it is not satisfied we can always add a constant to $\mathcal{E}_1$ without altering the original gradient flow.}
{Then} we introduce the so-called scalar auxiliary variable $r=\sqrt{\mathcal{E}_1}$,  
and reformulate the gradient flow \eqref{Gflow} as
\begin{subequations}\label{SAV1}
\begin{align}
 & \frac{\partial\phi}{\partial t}=\mathcal{G}\mu,\label{SAVeq_mu}\\
& \mu=\mathcal{L}\phi+\frac{r}{\sqrt{\mathcal{E}_1[\phi]}}U[\phi], \label{SAVeq_phi}\\
 & r_t=\frac{1}{2\sqrt{\mathcal{E}_1[\phi]}}
 \int_\Omega  U[\phi]\phi_t \, \mathrm{d}{\bf x}, \label{SAVeq_r}
\end{align}
\end{subequations}
where
\begin{equation}
  U[\phi]=\frac{\delta \mathcal{E}_1}{\delta \phi}. 
\end{equation}
{By} calculating the inner products of the above with $\mu$, $\frac{\partial\phi}{\partial t}$ and $2r$, respectively, we obtain the energy dissipation formulae for equation \eqref{SAV1}:
\begin{equation}
 \frac{d}{dt}[(\phi, \mathcal{L}\phi)+ r^2]=(\mu,\mathcal{G}\mu).
\end{equation}

\subsection{Evolutionary Kolmogorov–Arnold Networks}
The Evolutionary Kolmogorov–Arnold Networks (EvoKAN) follows the same structure as traditional evolutionary deep neural networks (EDNNs) \cite{du2021evolutional,zhang2024energy,zhang2022neural}. In particular, EvoKAN models Kolmogorov–Arnold Network weights as time-dependent
functions and updates them through the evolution of the governing
PDEs. Figure~\ref{fig:illustration} provides an overview of the EvoKAN. To explain EvoKAN, we start with the general nonlinear PDE~\eqref{eqn:general-pde}, where the solution approximation, $\bf u$, is represented using KANs with $L+1$ number of layers.
\begin{equation}
\hat{\boldsymbol{u}}(\boldsymbol{x}, \mathcal{K}(t)) = (\hat{u}_1, \hat{u}_2, \dots, \hat{u}_m)= KAN (\boldsymbol{x}, \mathcal{K}(t)).
\end{equation}
where $KAN (\boldsymbol{x}, \mathcal{K}(t))$ denotes $L+1$-layer KAN and  $\mathcal{K}(t)$, which varies with time, is the vector containing all parameters in the KAN and the input to the network is the spatial coordinate of the PDE:
\begin{equation}
\boldsymbol{x} = (x_1, x_2, \dots, x_d).   
\end{equation}
A fully connected KAN is defined in equation~\eqref{eq:kan}. Then the time evolution of the solution can be expressed as:
\begin{equation}
\frac{\partial \hat{\boldsymbol{u}}}{\partial t} = \frac{\partial \hat{\boldsymbol{u}}}{\partial \mathcal{K}} \frac{\partial \mathcal{K}}{\partial t}.
\end{equation}
At each time step, the time derivative of the network parameters, $\frac{\partial \mathcal{K}}{\partial t}$, is obtained by solving the optimization problem:
\begin{equation}\label{Wmin}
\frac{\partial \mathcal{K}}{\partial t} = argmin_{\gamma} \mathcal{J}(\gamma), \quad \mathcal{J}(\gamma) = \frac{1}{2} \int_{\Omega} \left\Vert \frac{\partial \hat{\boldsymbol{u}}}{\partial \mathcal{K}} \gamma - \mathcal{N}(\hat{\boldsymbol{u}})\right\Vert_2^2 \mathrm{d}\boldsymbol{x}.
\end{equation}
The optimality condition for this problem is:
\begin{equation}\label{OptCond}
\nabla_{\gamma} \mathcal{J}(\gamma_{\text{opt}}) = \left( \int_{\Omega} \frac{\partial \hat{\boldsymbol{u}}}{\partial \mathcal{K}}^{T} \frac{\partial \hat{\boldsymbol{u}}}{\partial \mathcal{K}} \mathrm{d}\boldsymbol{x} \right)\gamma_{\text{opt}} - \left( \int_{\Omega} \frac{\partial \hat{\boldsymbol{u}}}{\partial \mathcal{K}}^{T} \mathcal{N}(\hat{\boldsymbol{u}}) \mathrm{d}\boldsymbol{x} \right) = 0.
\end{equation}
The optimal $\gamma_{\text{opt}}$ is approximated by:
\begin{equation}\label{OptCondApprox}
\mathbf{J}^{T} \mathbf{J} \hat{\gamma}_{\text{opt}} = \mathbf{J}^{T} \mathbf{N},
\end{equation}
where $\mathbf{J}$ is the network's Jacobian and $\mathbf{N}$ is the PDE residual evaluated at collocation points:
\begin{equation}
\left( \mathbf{J} \right)_{ij} = \frac{\partial \hat{\boldsymbol{u}}^i}{\partial \mathcal{K}_j}, \quad \left( \mathbf{N} \right)_i = \mathcal{N}(\hat{\boldsymbol{u}}^i),
\end{equation}
where $i = 1, 2, \dots, N_u$ indexes the collocation points, and $j = 1, 2, \dots, N_{\mathcal{K}}$ indexes the network parameters. Automatic differentiation is used to compute the elements of $\mathbf{J}$ and $\mathbf{N}$. 
As a result, we can update the weights, $\mathcal{K}^{n+1}$ and the solution, $\hat{\boldsymbol{u}}^{n+1}(\boldsymbol{x}, \mathcal{K}(t))$ at $n+1$ time step as:
\begin{align}
\mathcal{K}^{n+1} &=\mathcal{K}^n +\hat{\gamma}_{\text{opt}}  \Delta t, \\
\hat{\boldsymbol{u}}^{n+1}(\boldsymbol{x}, \mathcal{K}(t))&=\hat{\boldsymbol{u}}(\boldsymbol{x}, \mathcal{K}^{n+1}(t)).
\end{align}

\section{Numerical Tests \label{sec:test}}
In what follows, we will present the numerical results for three different test cases: (1). one-dimensional Allen--Cahn equation
(2). two-dimensional Allen--Cahn equation and
(3). two--dimensional Navier--Stokes equations. 



\subsection{One-dimensional Allen--Cahn Equation}

In this section, we demonstrate the accuracy and performance of EvoKAN by solving the one-dimensional Allen--Cahn equation with periodic boundary conditions. The Allen--Cahn equation is a classical phase-field model used to describe phase separation processes in materials science \cite{allen1972ground,allen1973correction,shen2010numerical,feng2003numerical}. It can be viewed as a gradient flow of the Ginzburg--Landau free energy. In one-dimension case, the governing equation with boundary and initial conditions yields:
\begin{align}
    &u_t \;=\; \epsilon^2\,\frac{\partial^2 u}{\partial x^2} \;-\; g(u), \\
    &u(x,0) \;=\; a \,\sin\,\bigl(\pi x\bigr),\\
        &u(-L,t) \;=\;u(L,t), \, \forall t\ge 0
\end{align}
where $x\in\Omega = [-L,L]$ is the spatial coordinate, $t\ge 0$ is the time, and $\epsilon > 0$ is a small parameter controlling the thickness of the transition layer, and the nonlinearity is
\begin{align}
      g(u) \;=\; G'(u) \;=\; \frac{1}{\epsilon^2}\,u\,(u^2 - 1).
\end{align}
Here,
\begin{align}
   G(u) \;=\; \frac{1}{4\epsilon^2}\,(u^2 - 1)^2   
\end{align}
is the double-well potential that leads $u \approx \pm 1$.
In particular, the Allen--Cahn equation is the $L^2$-gradient flow of the Ginzburg--Landau free energy:
\begin{equation}
  E[u] \;=\; \int_\Omega \tfrac12\,\bigl|u_x\bigr|^2 \,dx 
  \;+\; \int_\Omega G(u)\,dx,
  \label{eq:ginzburg-landau-energy}
\end{equation}
where $G(u) = \tfrac{1}{4\epsilon^2}(u^2 - 1)^2$. When $\epsilon$ is small, $G(u)$ has two minima in $u=\pm1$, and the solution $u$ tends to form narrow transition layers (interfaces) separating these two stable states \cite{allen1972ground,allen1973correction}.
To evaluate our method, we used a high-accuracy spectral method as a benchmark \cite{shen1999wind,shen2010numerical}. We then compute a time-averaged $L^2$ error between the benchmark solution $u$ and our proposed solution $\widetilde{u}$ (from EvoKAN). Let $\Omega = [-1,1]$, and let $T$ be the final time. The error is defined as
\begin{equation}
  \mathcal{E}(u) 
  \;=\; \frac{1}{T}\,\int_{0}^T 
  \left( \frac{1}{|\Omega|}\,\int_\Omega 
    \bigl(u(x,t) - \widetilde{u}(x,t)\bigr)^2 \,dx
  \right)^{\!\frac12} dt.
  \label{eq:Allen--Cahn-error}
\end{equation}
\begin{figure}[H]
    \centering
    \begin{subfigure}[b]{0.32\linewidth}
        \centering
        \includegraphics[width=\linewidth]{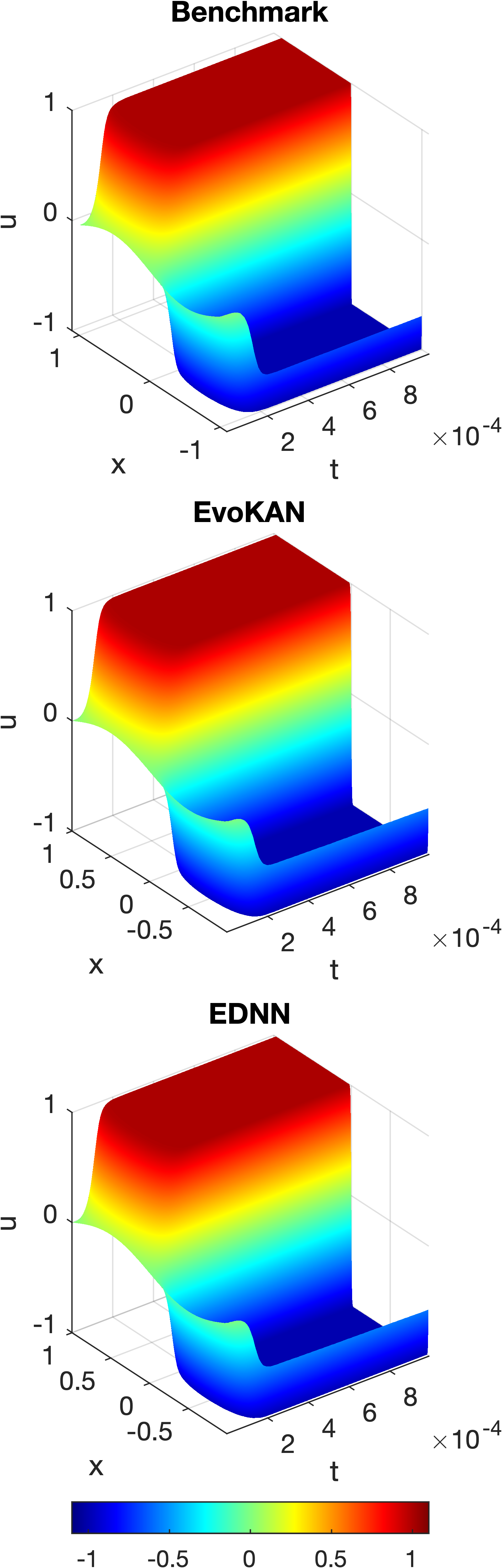}
        \caption{$\epsilon = 0.005$}
        \label{fig:ac-1d-eps0.005}
    \end{subfigure}
    \begin{subfigure}[b]{0.32\linewidth}
        \centering
        \includegraphics[width=\linewidth]{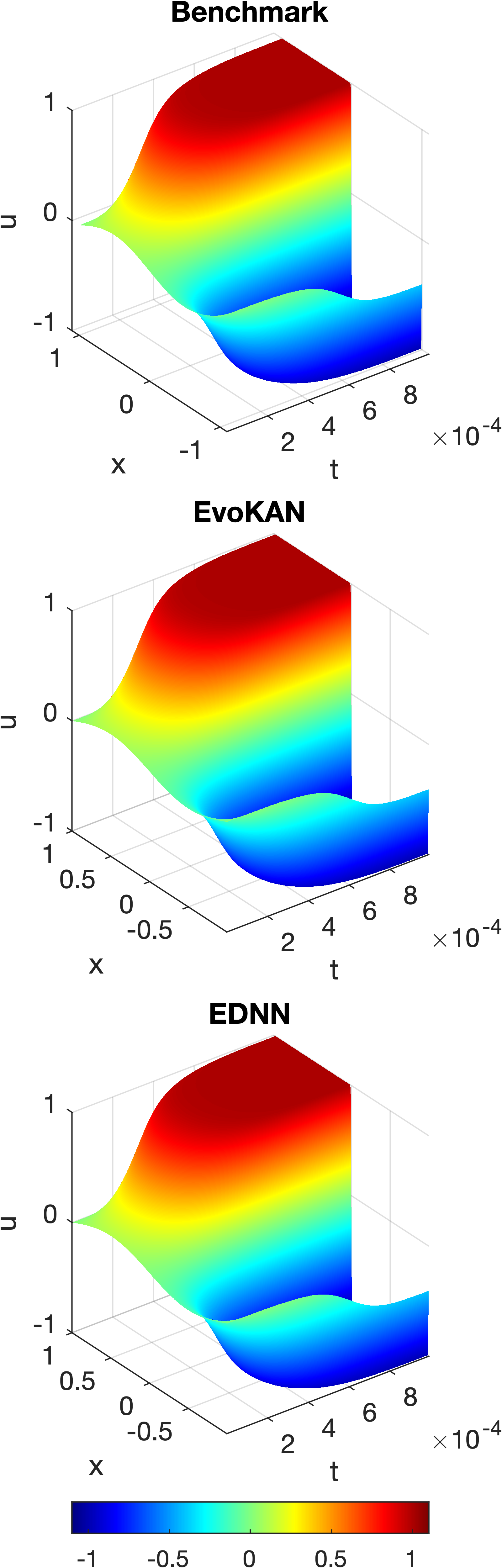}
        \caption{$\epsilon = 0.01$}
        \label{fig:ac-1d-eps0.01}
    \end{subfigure}
        \begin{subfigure}[b]{0.29\linewidth}
        \centering
        \includegraphics[width=\linewidth]{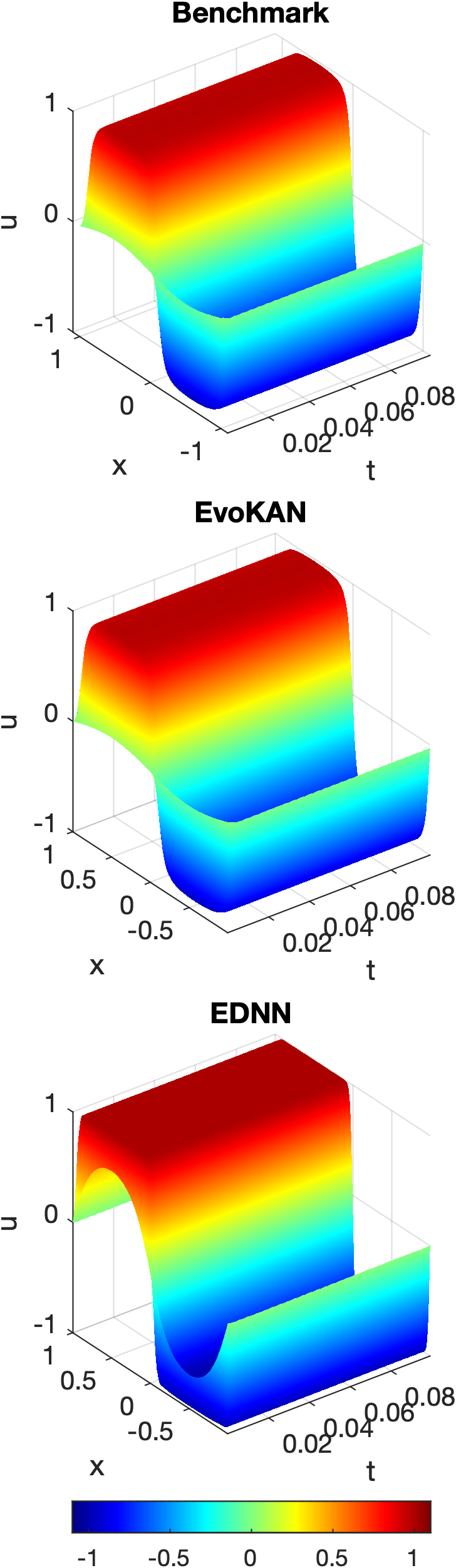}
        \caption{$\epsilon = 0.02$}
        \label{fig:ac-1d-eps0.02}
    \end{subfigure}
    \caption{
    Comparative analysis of benchmark solution and EvoKAN solution for one-dimensional Allen–Cahn equation. The three panels display the benchmark solution (top), EvoKAN solution (middle), and EDNN solution (bottom). 
    }
    \label{fig:ac-1d-1}
\end{figure}

We evaluate this error for multiple values of the parameter~$\epsilon$.
Figure \eqref{fig:ac-1d-1} shows the one-dimensional Allen--Cahn equation with the benchmark and the EvoKAN for different $\epsilon$ values.
Table~\ref{tab:Allen--Cahn-1d} shows the $L^2$ errors for the one-dimensional Allen--Cahn equation using our EvoKAN approach. The errors remain on the order of $10^{-4}$, demonstrating excellent agreement with the spectral benchmark across different scales of $\epsilon$.
\begin{table}[h]
    \centering
    \begin{tabular}{c|c|c|c}
      \hline
      Parameter
      & $\epsilon = 0.02$ 
      & $\epsilon = 0.01$
      & $\epsilon = 0.005$ 
       \\
      \hline
    EvoKAN $\mathcal{E}(u)$    
      & $1.6231\times 10^{-4}$ 
      & $2.3054\times 10^{-4}$       
      & $1.6236\times 10^{-4}$ \\
      \hline
    EDNN $\mathcal{E}(u)$    
      & $5.8736\times 10^{-4}$ 
      & $6.3897\times 10^{-4}$
      & $7.7592\times 10^{-4}$ 
       \\
      \hline      
    \end{tabular}
    \caption{$L^2$ errors (\eqref{eq:Allen--Cahn-error}) of EvoKAN and EDNN for the one-dimensional Allen--Cahn equation with various values of~$\epsilon$.}
    \label{tab:Allen--Cahn-1d}
\end{table}
As $\epsilon$ decreases, the interface becomes sharper and the solution transitions more abruptly between the stable phases $u \approx \pm 1$. 
This can pose numerical challenges because a finer spatial resolution is required to capture steep gradients accurately. Nonetheless, the EvoKAN maintains a stable error level for these increasingly stiff problems, indicating its robustness in handling sharp interfaces.

The Allen--Cahn equation serves as a foundational test problem in phase-field modeling, where capturing the equilibria governed by the Ginzburg--Landau free energy is essential. The close agreement with the spectral method illustrates that EvoKAN can faithfully follow the gradient-flow dynamics driving the system to minimize the Ginzburg--Landau free energy (cf.~\eqref{eq:ginzburg-landau-energy}). Consequently, the computed solutions accurately reflect the expected behavior of forming and moving interfaces, ultimately settling into near-equilibrium states as governed by the potential $G(u)$.
Overall, these results show that the proposed EvoKAN is effective to solve the one-dimensional phase-field problems.

%



\subsection{Two-dimensional Allen--Cahn Equation}

The two-dimensional Allen--Cahn equation extends the phase-field description to a two-dimensional spatial domain, where the order parameter $u(\mathbf{x},t)$ now depends on $\mathbf{x}=(x,y)\in \Omega = [-L,L]\times[-L,L]$. The governing equation with initial and boundary conditions is
\begin{align}
    &u_t \;=\; \epsilon^2\,\Delta u \;-\; g(u),
    \\
    &u(x,y,0) \;=\; 0.08 \,\sin\,\bigl(\alpha\pi x\bigr)\,\sin\,\bigl(\alpha \pi y\bigr),\label{eqn:ic-2d-ac}
    \\
    & u(L,L,t)  \;=\; u(L,-L,t),   u(L,L,t)  \;=\; u(-L,L,t),\quad \forall t\ge 0
\end{align}
where $\Delta$ denotes the Laplace operator in two spatial dimensions,  $\alpha$ is an integer parameter chosen to guarantee a non-trivial interface evolution, and $g(u)$ is the derivative of the double-well potential $G(u) = \tfrac{1}{4\epsilon^2}(u^2-1)^2$. As in the one-dimensional case, the solution $u$ evolves to minimize the Ginzburg--Landau free energy, but the resulting interfaces between the stable phases $u \approx \pm 1$ now appear as curves (or contours) in the plane \cite{allen1976mechanisms,sun2007sharp,chiu2011conservative}. These interfaces exhibit curvature-driven motion, tending to shorten their overall length to reduce the free energy. Consequently, numerical simulations in two dimensions typically involve interface shrinkage, collision, and annihilation, leading to coarsening of the domain morphology. 
Due to the steep gradients within the transition regions, accurate resolution of the two-dimensional Allen--Cahn equation requires careful spatial discretization. 
As such, it remains a fundamental testbed for phase-field methods aiming to capture complex interface dynamics in higher-dimensional settings.
In this section, we consider the Allen--Cahn equation in a two-dimensional square domain $(x,y \in [-1,1]$ and fix a time step size $\Delta t$. Periodic boundary conditions are imposed in both spatial directions. The parameter in the double-well potential is set to $\epsilon = 0.05$.

\begin{figure}[H]
    \centering
    \begin{subfigure}[b]{\linewidth}
        \centering
        \includegraphics[width=\linewidth]{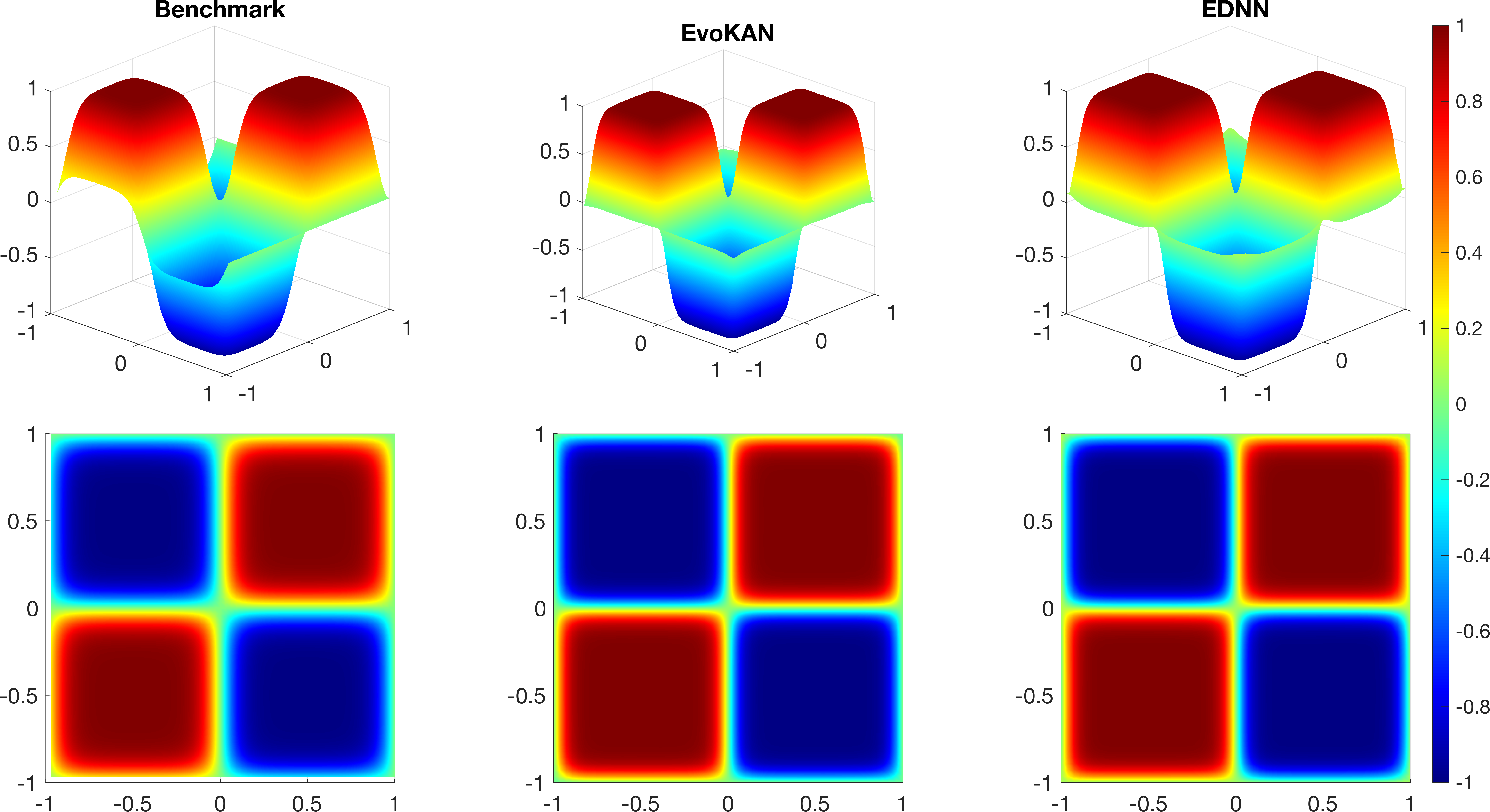}
        \caption{$\epsilon = 0.05$, initial condition \eqref{eqn:ic-2d-ac} with $\alpha=1$.}
        \label{fig:ac-2d-ic1}
    \end{subfigure}
    \begin{subfigure}[b]{\linewidth}
        \centering
        \includegraphics[width=\linewidth]{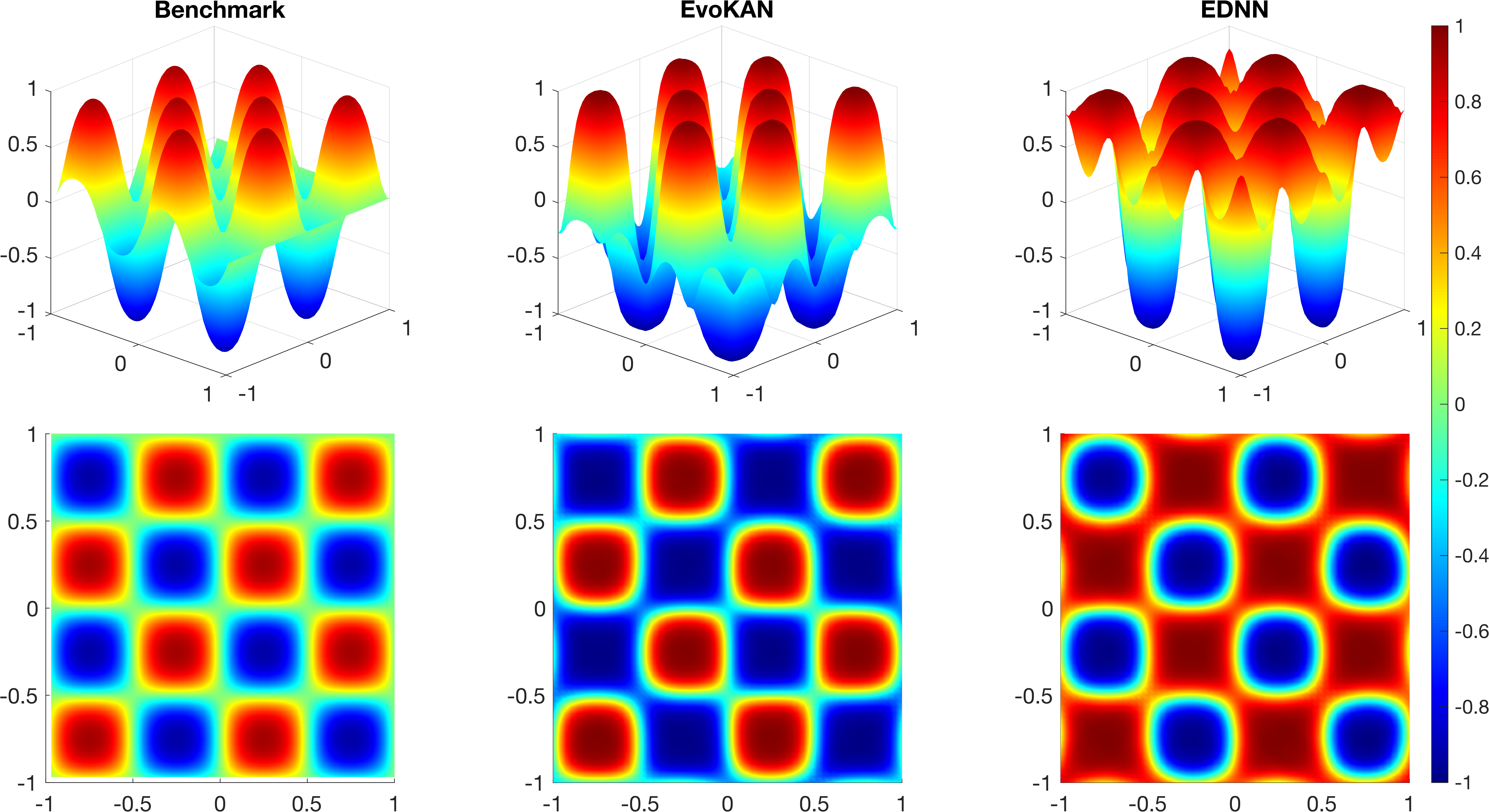}
        \caption{$\epsilon = 0.05$, initial condition \eqref{eqn:ic-2d-ac} with $\alpha=2$.}
        \label{fig:ac-2d-ic2}
    \end{subfigure}
    \caption{
    Comparative analysis of benchmark solution, EvoKAN, and EDNN solutions for the two-dimensional Allen–Cahn equation under different initial conditions. Panels (a) and (b) show the results for initial conditions with $\alpha = 1$ and $\alpha = 2$, respectively. 
    }
    \label{fig:ac-2d-1}
\end{figure}

Figure~\ref{fig:ac-2d-1} presents the numerical solutions of the two-dimensional Allen--Cahn equation under various initial conditions, comparing the benchmark solution with the EvoKAN and EDNN approaches. Across all tested scenarios, the EvoKAN solutions remain in close agreement with the benchmark, effectively capturing the curvature-driven interface dynamics and coarsening phenomena characteristic of Allen--Cahn evolution. In particular, the spatial distribution of phases and the location of transition regions are well-matched between the EvoKAN and the benchmark. The EDNN solutions also yield plausible phase-field structures, though minor deviations in interface shape or position may appear in select cases. Overall, EvoKAN is robust in reproducing the benchmark solution for different initial conditions, indicating its strong potential for accurately capturing phase-transition behavior in two-dimensional domains.




\subsection{Navier-Stokes Equation \label{ss-nse}}
The incompressible fluid are governed by the Navier-Stokes equations, which describe the conservation of momentum and mass in a fluid system \cite{constantin1988navier,temam2024navier,chorin1968numerical}. { The governing equations of NSE in $\Omega\times [0,T]$ yield the following \cite{layton2008introduction}:
\begin{align} 
    &\frac{\partial \mathbf{v}}{\partial t} + (\mathbf{v} \cdot \nabla)\mathbf{v} = -\nabla p + \nu \nabla^2 \mathbf{v} + \mathbf{f},\label{eq:navier_stokes}
    \\
    &\nabla \cdot \mathbf{v} = 0,\label{eq:incompressibility}
    \\
     &u(x,y,0) = - \sin(2\pi y), v(x,y,0) = \cos(2\pi y), \label{nse-ics}\\
         & u(L,L,t)  \;=\; u(L,-L,t),   u(L,L,t)  \;=\; u(-L,L,t),\label{nse-bc-1}
         \\
             & v(L,L,t)  \;=\; v(L,-L,t),   v(L,L,t)  \;=\; v(-L,L,t),\quad \forall t\ge 0.\label{nse-bc-2}
\end{align}
Among them, equation \eqref{eq:navier_stokes} is
the momentum equation, which captures the time evolution of the velocity field, is expressed as
where $\mathbf{v} = (u, v)$ represents the velocity field of the fluid, with $u$ and $v$ denoting the velocity components in the $x$- and $y$-directions, respectively. The term $p$ is the pressure field, which enforces the incompressibility constraint by redistributing forces within the fluid. The parameter $\nu$ is the kinematic viscosity of the fluid, which characterizes the resistance to shear and momentum diffusion.
The incompressibility condition for the NSE is imposed in equation ~\eqref{eq:incompressibility},
which ensures that the fluid's density remains constant. Equation \eqref{nse-ics} imposes the initial conditions for velocities and equations \eqref{nse-bc-1}--\eqref{nse-bc-2} indicates the doubly periodic boundary conditions. 
\begin{figure}[H]
    \centering
    \includegraphics[width=\linewidth]{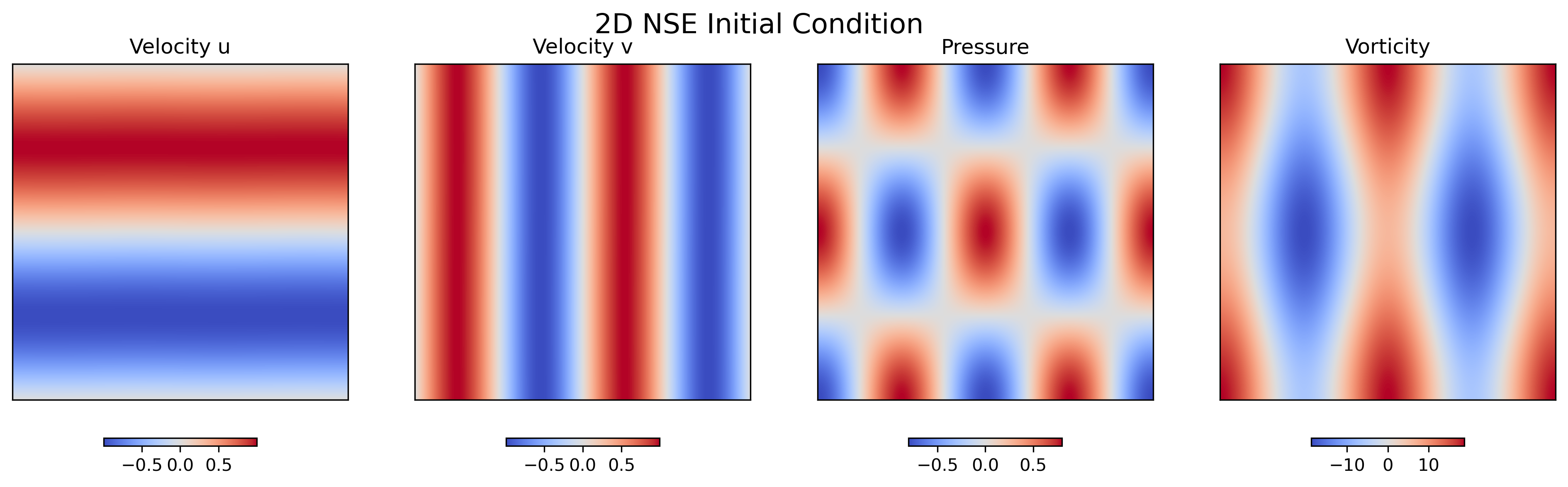}
    \caption{
     Initial conditions for the two-dimensional Navier-Stokes equations. The figure presents the initial velocity fields ($u$) (first column) and ($v$) (second column), pressure field ($P$) (third column), and vorticity ($\omega$) (fourth column).}
    \label{fig:NSE-nu-ic}
\end{figure}
In the numerical tests, we consider the square domain $\Omega = [-1,1]\times[-1,1]$ with double periodic boundary conditions along with the following initial conditions.}
For simplicity, the forcing term $\mathbf{f}$ is set to zero. 
We use the spectral method to generate full-order model solutions as a benchmark. The details of the implementation are given in Section~\ref{appendix-nse-spectral}.
Figure~\ref{fig:NSE-nu-ic} shows the initial conditions for the velocity and pressure fields as well as the vorticity. 

We consider two different test regimes: (1). $\nu=0.05$ and (2). $\nu=0.01$.

\begin{figure}[H]
    \begin{subfigure}[b]{\linewidth}
        \centering
    \includegraphics[width=\linewidth]{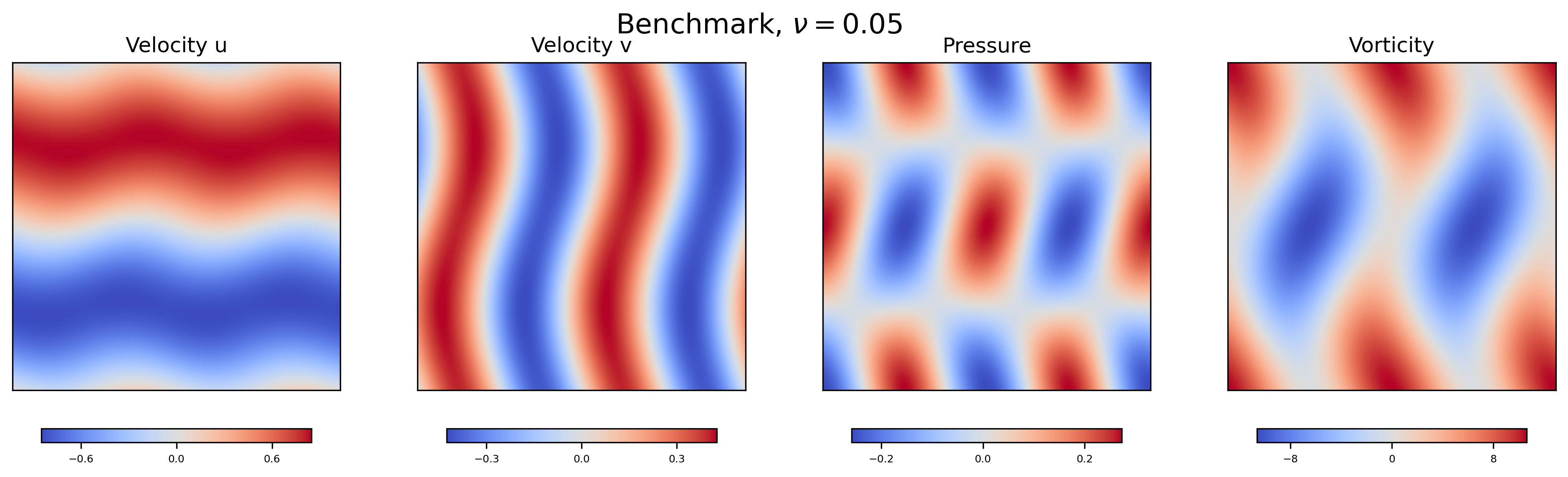}
\includegraphics[width=\linewidth]{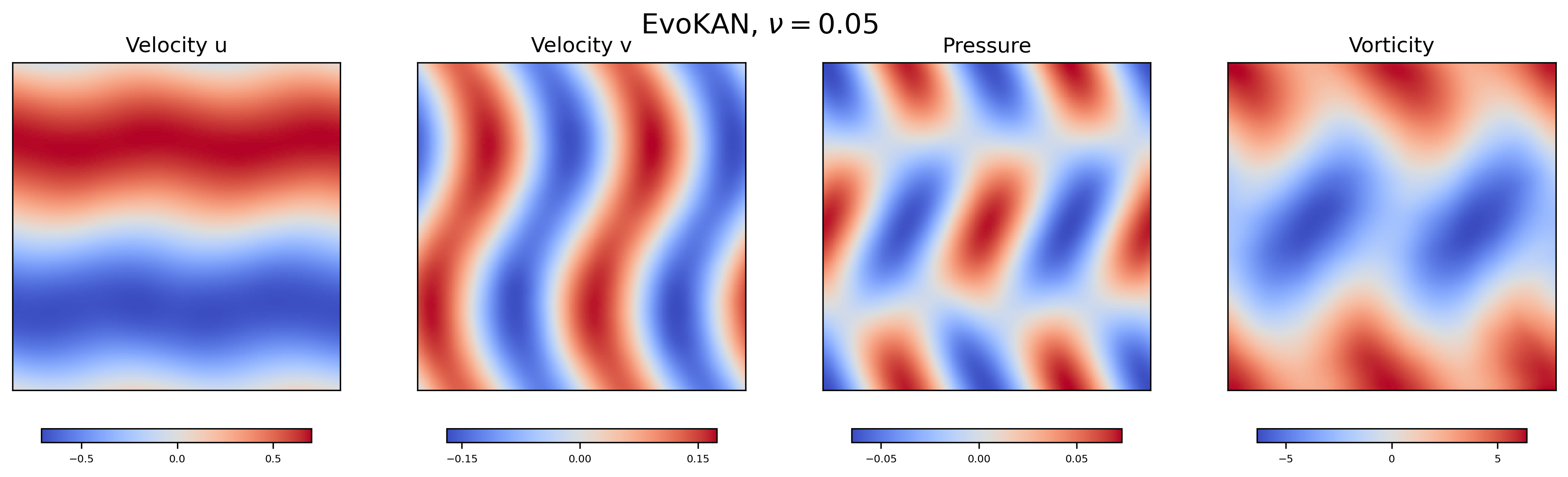}

\caption{$t=0.1$, (top) Benchmark, (bottom) EvoKAN result}
        \label{fig:nse-2d-nu1-t1}
    \end{subfigure}    
     \begin{subfigure}[b]{\linewidth}
        \centering
    \includegraphics[width=\linewidth]{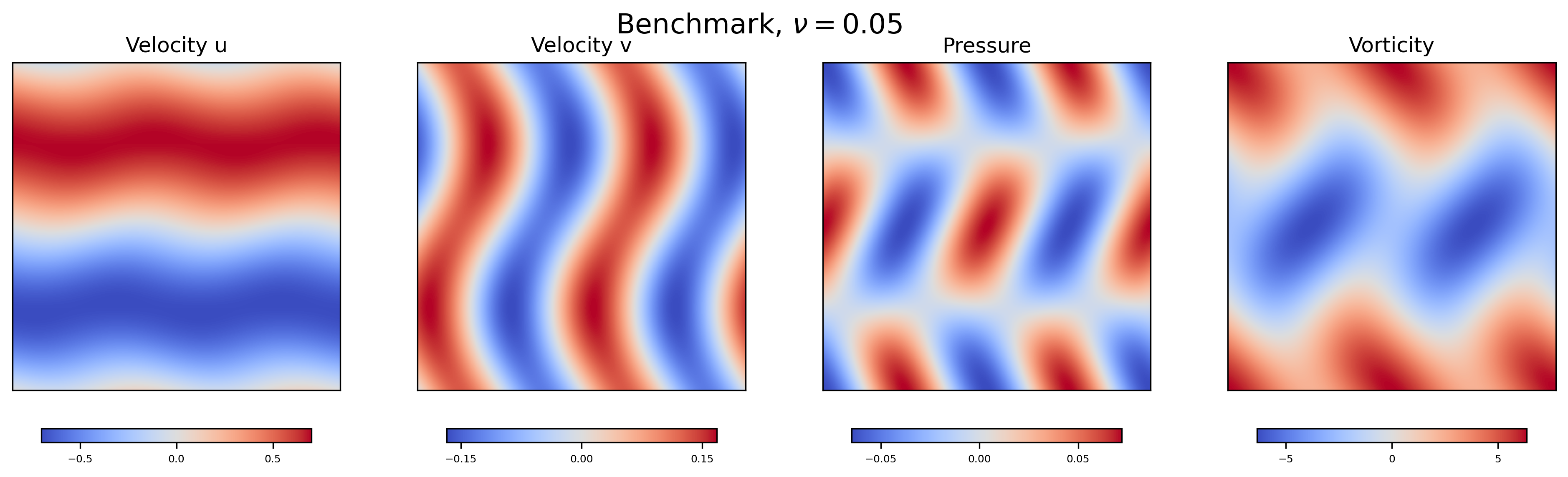}
    \includegraphics[width=\linewidth]{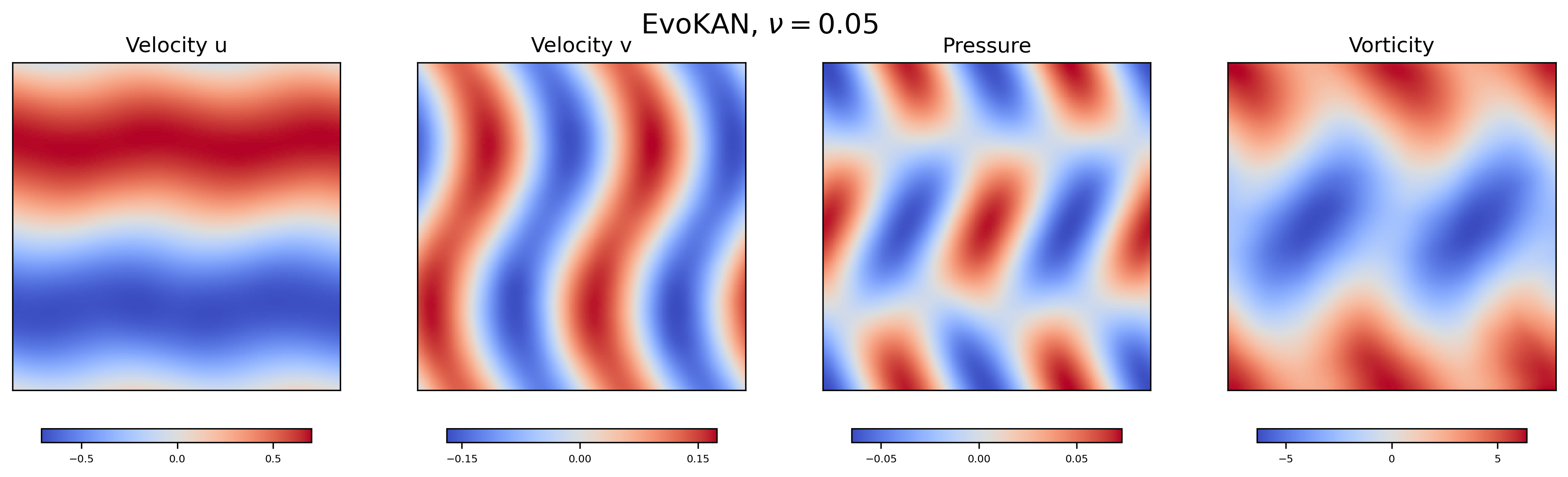}
\caption{$t=0.2$, (top) Benchmark, (bottom) EvoKAN result}
        \label{fig:nse-2d-nu1-t2}
    \end{subfigure}      
    \caption{Comparison of benchmark and EvoKAN for two-dimension NSE with $\nu=0.05$: Velocity fields, $u$ (first column) and $v$ (second column) , pressure $P$ (third column) , and vorticity $\omega$ (fourth column).}
    \label{fig:NSE-nu-1}
\end{figure}

\begin{figure}[H]
    \begin{subfigure}[b]{\linewidth}
        \centering
\includegraphics[width=\linewidth]{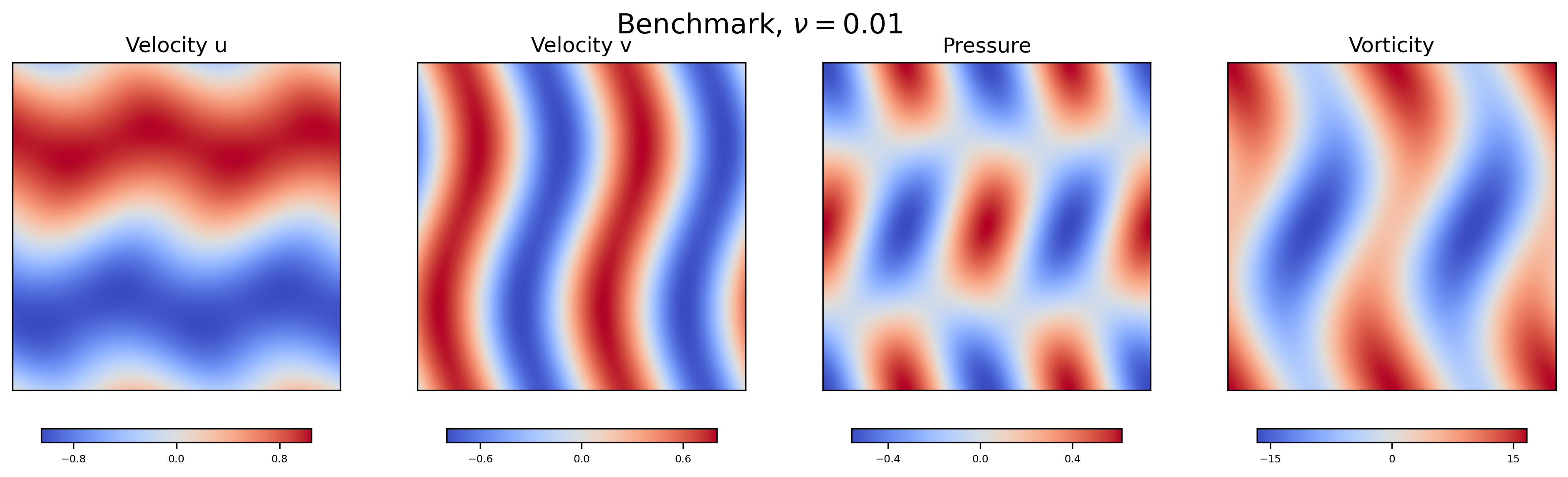}
\includegraphics[width=\linewidth]{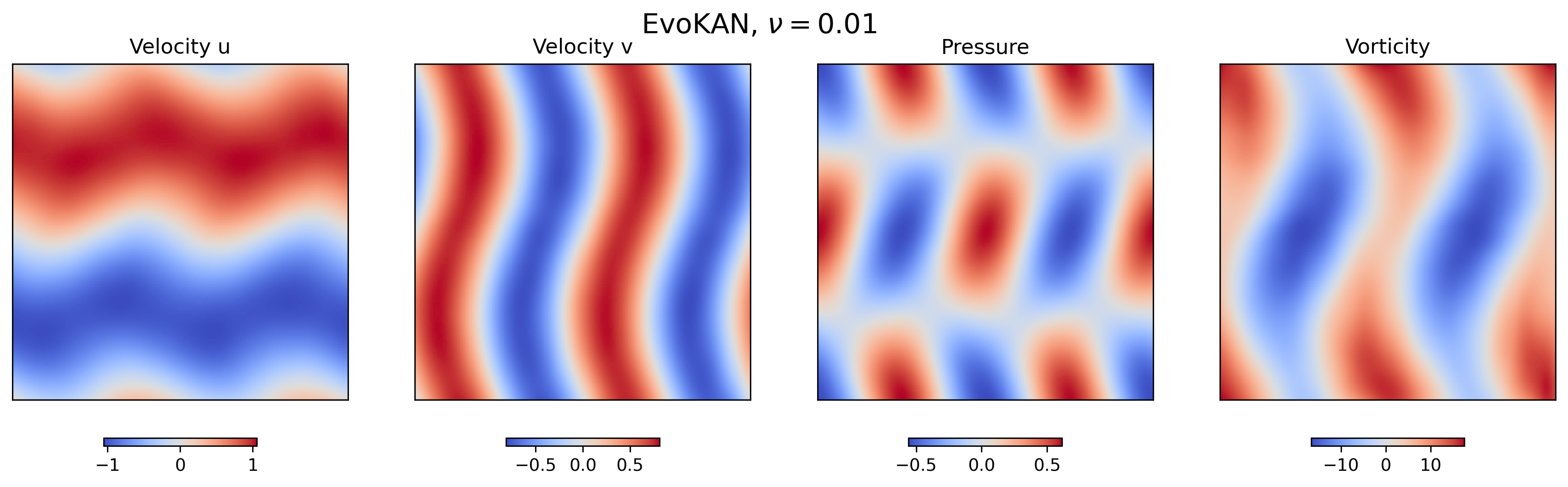}
\caption{$t=0.1$, (top) Benchmark, (bottom) EvoKAN result}
        \label{fig:nse-2d-nu2-t1}
    \end{subfigure}    
     \begin{subfigure}[b]{\linewidth}
        \centering
    \includegraphics[width=\linewidth]{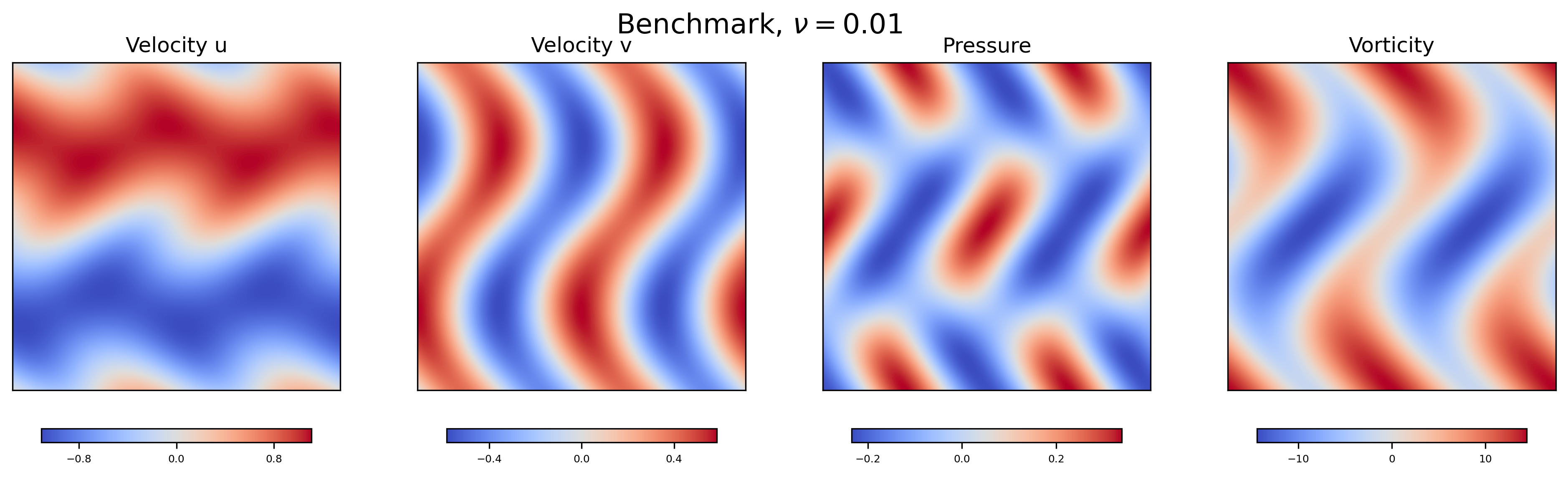}
    \includegraphics[width=\linewidth]{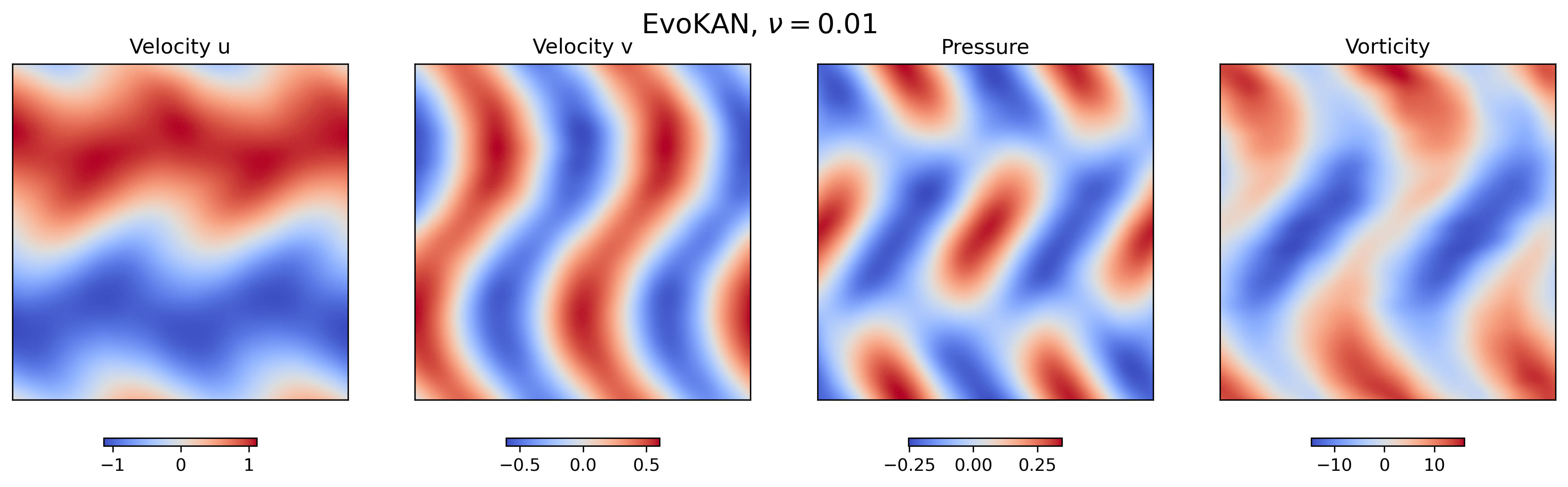}
\caption{$t=0.2$, (top) Benchmark, (bottom) EvoKAN result}
        \label{fig:nse-2d-nu2-t2}
    \end{subfigure}      
    \caption{Comparison of benchmark and EvoKAN for two-dimension NSE with $\nu=0.01$: Velocity fields, $u$ (first column) and $v$ (second column) , pressure $P$ (third column) , and vorticity $\omega$ (fourth column).}
    \label{fig:NSE-nu-2}
\end{figure}

The EvoKAN model demonstrates a strong capacity to capture and reproduce the underlying flow structures of the incompressible Navier-Stokes equations, as evidenced in Figures~\ref{fig:NSE-nu-1} and~\ref{fig:NSE-nu-2}. In these figures, we compare the velocity fields ($u$ and $v$), pressure $P$, and vorticity $\omega$ for EvoKAN and the benchmark solutions generated by the spectral method. For both viscosity values considered ($\nu=0.05$ and $\nu=0.01$), EvoKAN shows good agreement with the benchmark solution, demonstrating its ability to approximate the transient dynamics of the two-dimensional incompressible flow. It is worth noting that EDNN fails to reproduce the unstable solutions for both regimes that only EvoKAN and benchmark solutions are compared. 

In particular, in the case of moderate viscosity ($\nu=0.05$, Figure~\ref{fig:NSE-nu-1}), the flow remains relatively smooth, which helps to stabilize the flow structures. Here, EvoKAN accurately follows the benchmark solution, reproducing the dominant vortex patterns, flow rotations, and pressure distributions. The vorticity contours in particular align well with the benchmark results, indicating that the EvoKAN retains key nonlinear advective and diffusive processes inherent in the governing Navier-Stokes equations.
As the viscosity decreases to $\nu=0.01$ (Figure~\ref{fig:NSE-nu-2}), the flow becomes more energetic and develops sharper gradients, making it more challenging to approximate. However, EvoKAN continues to reproduce the main flow features and the overall vorticity evolution.
In summary, the numerical results show that EvoKAN is effective in solving complex fluid problems.

\section{Conclusion and Future Work\label{sec:conclusion}}
This paper introduces a novel Evolutionary Kolmogorov–Arnold Network to solve time-dependent partial differential equations. Unlike traditional approaches that rely on multilayer perceptrons (MLPs), EvoKAN employs spline-based basis functions \cite{liu2024kan,liu2024kan2}, offering enhanced flexibility, smoothness, and localized control in representing complex spatial domains. Automatic differentiation is used to compute spatial derivatives directly from the spline-based network output.

EvoKAN has the following four significant advantages. (1) \textbf{Simple network and no need for global spatio-temporal optimization:} Traditional neural network-based methods for time-dependent PDE often require global optimization throughout the spatio-temporal domain \cite{meng2020ppinn,chen2023implicit}. In contrast, EvoKAN’s state represents only the instantaneous PDE solution, significantly reducing network complexity for a given problem. (2) \textbf{Explicit time dependence and causality:} EvoKAN the Kolmogorov–Arnold network weights as time-dependent functions and updates them through the evolution of the governing
PDEs. EvoKAN circumvents the limitations of methods that only minimize equation residuals. (3) \textbf{Long-horizon prediction:} Crucially, EvoKAN can predict very long time horizons, even in chaotic regimes, overcoming a major challenge that has hindered other neural network-based PDE solvers \cite{lai2018modeling,long2019pde,blechschmidt2021three}.
(4) \textbf{Exceptional versatility and accuracy:} EvoKAN’s versatility and accuracy have been demonstrated on several PDE problems, including the one-dimensional and two-dimensional Allen–Cahn equations, and the two-dimensional incompressible Navier–Stokes equations. In all cases, the accuracy of EvoKAN improved monotonically as the structure of the network was refined and the spatio-temporal resolution increased.

Building on EvoKAN's effectiveness in capturing complex dynamics and accommodating spatiotemporal flexibility, several future directions are necessary for further exploration. 
First, an adaptive strategy that automatically refines the spline-based network architecture in regions of steep gradients or highly localized features could improve both accuracy and efficiency, especially for problems with sharp transitions or multiscale phenomena \cite{wong2010adaptive,he2015adaptive,ge1999adaptive,mou2021data,harlim2021machine,liang2022stiffness,popov2021multifidelity}. Second, extending EvoKAN to more complex geometries and boundary conditions through domain decomposition or advanced spline-based representations would broaden its applicability to real-world engineering and scientific problems \cite{berg2018unified,winovich2019convpde}. A third direction involves the integration of reduced order model ROMs, data assimilation, or state estimation methodologies, enabling EvoKAN to incorporate noisy or partial observations in real time and to tackle uncertainty quantification problems \cite{deshmukh2016neural,mou2023combining,mou2023efficient}. Finally, coupling EvoKAN with high-performance computing and parallelization strategies could greatly accelerate large-scale simulations, shedding light on real-time PDE solvers for high-resolution simulations in climate science and engineering applications \cite{pijanowski2014big,yan2021accelerating}.

\section*{Acknowledgment}
We would like to thank the support of National Science Foundation (DMS-2053746, DMS-2134209, ECCS-2328241, CBET-2347401 and OAC-2311848), and U.S.~Department of Energy (DOE) Office of Science Advanced Scientific Computing Research program DE-SC0023161, and DOE–Fusion Energy Science, under grant number: DE-SC0024583.

\bibliographystyle{plainnat} 
\bibliography{ref,ref_qg,ref_kan}

\appendix

\section{Spectral Method for 2D Incompressible NSE \label{appendix-nse-spectral}}
In what follows, we briefly introduce the spectral method for solving the two-dimensional Navier-Stokes equations used in Section~\ref{ss-nse}.
\subsection{Momentum Equation}
The Navier-Stokes equations for an incompressible fluid are given by:
\begin{equation}
    \frac{\partial \mathbf{v}}{\partial t} + (\mathbf{v} \cdot \nabla)\mathbf{v} = -\nabla p + \nu \nabla^2 \mathbf{v} + \mathbf{f},
\end{equation}
where:
\begin{itemize}
    \item \(\mathbf{v} = (v_x, v_y)\): Velocity field.
    \item \(p\): Pressure field.
    \item \(\nu\): Kinematic viscosity of the fluid.
    \item \(\nabla^2\): Laplace operator, representing diffusion.
    \item \(\mathbf{f}\): External forcing term (not considered in this example).
\end{itemize}

\subsection{Incompressibility Condition (Mass Conservation)}
The incompressibility condition ensures the fluid density remains constant:
\begin{equation}
    \nabla \cdot \mathbf{v} = 0.
\end{equation}

In Fourier space, this condition becomes:
\begin{equation}
    i\mathbf{k} \cdot \hat{\mathbf{v}} = 0,
\end{equation}
where \(\mathbf{k}\) is the wave vector, enforcing orthogonality between the velocity field and the wave vector.

\subsection{Fourier Transform}
For a periodic velocity field \(\mathbf{v}(x, y, t)\), the Fourier expansion is:
\begin{equation}
    \mathbf{v}(x, y, t) = \sum_{k_x, k_y} \hat{\mathbf{v}}(k_x, k_y, t) e^{i(k_x x + k_y y)},
\end{equation}
where:
\begin{itemize}
    \item \(\hat{\mathbf{v}}(k_x, k_y, t)\): Fourier coefficients of the velocity field.
    \item \(k_x, k_y\): Wavenumbers in the \(x\)- and \(y\)-directions.
\end{itemize}

Differential operators are simplified in Fourier space:
\begin{itemize}
    \item Gradient operator \(\nabla\): \(i\mathbf{k}\), where \(\mathbf{k} = (k_x, k_y)\).
    \item Laplace operator \(\nabla^2\): \(-|\mathbf{k}|^2\).
\end{itemize}

\subsection{Momentum Equation in Fourier Space}
The Navier-Stokes momentum equation in Fourier space is:
\begin{equation}
    \frac{\partial \hat{\mathbf{v}}}{\partial t} + \widehat{(\mathbf{v} \cdot \nabla)\mathbf{v}} = -i\mathbf{k} \hat{p} - \nu |\mathbf{k}|^2 \hat{\mathbf{v}}.
\end{equation}

The terms are represented as:
\begin{itemize}
    \item \textbf{Time Derivative:}
    \begin{equation}
        \frac{\partial \mathbf{v}}{\partial t} \xrightarrow{\text{Fourier Transform}} \frac{\partial \hat{\mathbf{v}}}{\partial t}.
    \end{equation}
    \item \textbf{Advection Term:}
    \begin{equation}
        (\mathbf{v} \cdot \nabla)\mathbf{v} \xrightarrow{\text{Fourier Transform}} \text{Convolution of Fourier modes}.
    \end{equation}
    \item \textbf{Pressure Gradient:}
    \begin{equation}
        -\nabla p \xrightarrow{\text{Fourier Transform}} -i\mathbf{k} \hat{p}.
    \end{equation}
    \item \textbf{Viscous Diffusion:}
    \begin{equation}
        \nu \nabla^2 \mathbf{v} \xrightarrow{\text{Fourier Transform}} -\nu |\mathbf{k}|^2 \hat{\mathbf{v}}.
    \end{equation}
\end{itemize}

\subsection{Pressure Elimination}
The incompressibility condition is used to eliminate the pressure term:
\begin{equation}
    \hat{p} = \frac{\mathbf{k} \cdot \widehat{(\mathbf{v} \cdot \nabla)\mathbf{v}}}{|\mathbf{k}|^2}.
\end{equation}

\subsection{Dealiasing}
Nonlinear terms, such as \((\mathbf{v} \cdot \nabla)\mathbf{v}\), introduce high-frequency noise due to mode convolution. To avoid aliasing, the \textbf{2/3-rule} is applied:
\begin{equation}
    |\mathbf{k}| < \frac{2}{3} k_{\text{max}},
\end{equation}
where high-frequency components are truncated for numerical stability.

%% file: arxiv_draft_submission.bbl
\begin{thebibliography}{69}
\providecommand{\natexlab}[1]{#1}
\providecommand{\url}[1]{\texttt{#1}}
\expandafter\ifx\csname urlstyle\endcsname\relax
  \providecommand{\doi}[1]{doi: #1}\else
  \providecommand{\doi}{doi: \begingroup \urlstyle{rm}\Url}\fi

\bibitem[Abueidda et~al.(2024)Abueidda, Pantidis, and
  Mobasher]{abueidda2024deepokan}
Diab~W Abueidda, Panos Pantidis, and Mostafa~E Mobasher.
\newblock Deepokan: Deep operator network based on kolmogorov arnold networks
  for mechanics problems.
\newblock \emph{arXiv preprint arXiv:2405.19143}, 2024.

\bibitem[Allen and Cahn(1973)]{allen1973correction}
Samuel~M Allen and John~W Cahn.
\newblock A correction to the ground state of fcc binary ordered alloys with
  first and second neighbor pairwise interactions.
\newblock \emph{Scripta Metallurgica}, 7\penalty0 (12):\penalty0 1261--1264,
  1973.

\bibitem[Allen and Cahn(1976)]{allen1976mechanisms}
Samuel~M Allen and John~W Cahn.
\newblock Mechanisms of phase transformations within the miscibility gap of
  fe-rich fe-al alloys.
\newblock \emph{Acta Metallurgica}, 24\penalty0 (5):\penalty0 425--437, 1976.

\bibitem[Allen and Cahn(1972)]{allen1972ground}
Samuel~Miller Allen and John~W Cahn.
\newblock Ground state structures in ordered binary alloys with second neighbor
  interactions.
\newblock \emph{Acta Metallurgica}, 20\penalty0 (3):\penalty0 423--433, 1972.

\bibitem[Apicella et~al.(2021)Apicella, Donnarumma, Isgr{\`o}, and
  Prevete]{apicella2021survey}
Andrea Apicella, Francesco Donnarumma, Francesco Isgr{\`o}, and Roberto
  Prevete.
\newblock A survey on modern trainable activation functions.
\newblock \emph{Neural Networks}, 138:\penalty0 14--32, 2021.

\bibitem[Beck and Kurz(2021)]{beck2021perspective}
Andrea Beck and Marius Kurz.
\newblock A perspective on machine learning methods in turbulence modeling.
\newblock \emph{GAMM-Mitteilungen}, 44\penalty0 (1):\penalty0 e202100002, 2021.

\bibitem[Berg and Nystr{\"o}m(2018)]{berg2018unified}
Jens Berg and Kaj Nystr{\"o}m.
\newblock A unified deep artificial neural network approach to partial
  differential equations in complex geometries.
\newblock \emph{Neurocomputing}, 317:\penalty0 28--41, 2018.

\bibitem[Blechschmidt and Ernst(2021)]{blechschmidt2021three}
Jan Blechschmidt and Oliver~G Ernst.
\newblock Three ways to solve partial differential equations with neural
  networks—a review.
\newblock \emph{GAMM-Mitteilungen}, 44\penalty0 (2):\penalty0 e202100006, 2021.

\bibitem[Boull{\'e} and Townsend(2023)]{boulle2023mathematical}
Nicolas Boull{\'e} and Alex Townsend.
\newblock A mathematical guide to operator learning.
\newblock \emph{arXiv preprint arXiv:2312.14688}, 2023.

\bibitem[Braun and Griebel(2009)]{braun2009constructive}
J{\"u}rgen Braun and Michael Griebel.
\newblock On a constructive proof of kolmogorov’s superposition theorem.
\newblock \emph{Constructive approximation}, 30:\penalty0 653--675, 2009.

\bibitem[Cai et~al.(2021{\natexlab{a}})Cai, Mao, Wang, Yin, and
  Karniadakis]{cai2021physics}
Shengze Cai, Zhiping Mao, Zhicheng Wang, Minglang Yin, and George~Em
  Karniadakis.
\newblock Physics-informed neural networks (pinns) for fluid mechanics: A
  review.
\newblock \emph{Acta Mechanica Sinica}, 37\penalty0 (12):\penalty0 1727--1738,
  2021{\natexlab{a}}.

\bibitem[Cai et~al.(2021{\natexlab{b}})Cai, Wang, Lu, Zaki, and
  Karniadakis]{cai2021deepm}
Shengze Cai, Zhicheng Wang, Lu~Lu, Tamer~A Zaki, and George~Em Karniadakis.
\newblock Deepm\&mnet: Inferring the electroconvection multiphysics fields
  based on operator approximation by neural networks.
\newblock \emph{Journal of Computational Physics}, 436:\penalty0 110296,
  2021{\natexlab{b}}.

\bibitem[Chen et~al.(2023)Chen, Wu, Grinspun, Zheng, and
  Chen]{chen2023implicit}
Honglin Chen, Rundi Wu, Eitan Grinspun, Changxi Zheng, and Peter~Yichen Chen.
\newblock Implicit neural spatial representations for time-dependent pdes.
\newblock In \emph{International Conference on Machine Learning}, pages
  5162--5177. PMLR, 2023.

\bibitem[Chiu and Lin(2011)]{chiu2011conservative}
Pao-Hsiung Chiu and Yan-Ting Lin.
\newblock A conservative phase field method for solving incompressible
  two-phase flows.
\newblock \emph{Journal of Computational Physics}, 230\penalty0 (1):\penalty0
  185--204, 2011.

\bibitem[Chorin(1968)]{chorin1968numerical}
Alexandre~Joel Chorin.
\newblock Numerical solution of the navier-stokes equations.
\newblock \emph{Mathematics of computation}, 22\penalty0 (104):\penalty0
  745--762, 1968.

\bibitem[Constantin and Foia{\c{s}}(1988)]{constantin1988navier}
Peter Constantin and Ciprian Foia{\c{s}}.
\newblock \emph{Navier-stokes equations}.
\newblock University of Chicago press, 1988.

\bibitem[Cuomo et~al.(2022)Cuomo, Di~Cola, Giampaolo, Rozza, Raissi, and
  Piccialli]{cuomo2022scientific}
Salvatore Cuomo, Vincenzo~Schiano Di~Cola, Fabio Giampaolo, Gianluigi Rozza,
  Maziar Raissi, and Francesco Piccialli.
\newblock Scientific machine learning through physics--informed neural
  networks: Where we are and what’s next.
\newblock \emph{Journal of Scientific Computing}, 92\penalty0 (3):\penalty0 88,
  2022.

\bibitem[Deshmukh et~al.(2016)Deshmukh, Deo, Bhaskaran, Nair, and
  Sandhya]{deshmukh2016neural}
Aditya~N Deshmukh, MC~Deo, Prasad~K Bhaskaran, TM~Balakrishnan Nair, and
  KG~Sandhya.
\newblock Neural-network-based data assimilation to improve numerical ocean
  wave forecast.
\newblock \emph{IEEE Journal of Oceanic Engineering}, 41\penalty0 (4):\penalty0
  944--953, 2016.

\bibitem[Du and Zaki(2021)]{du2021evolutional}
Yifan Du and Tamer~A Zaki.
\newblock Evolutional deep neural network.
\newblock \emph{Physical Review E}, 104\penalty0 (4):\penalty0 045303, 2021.

\bibitem[Fakhoury et~al.(2022)Fakhoury, Fakhoury, and
  Speleers]{fakhoury2022exsplinet}
Daniele Fakhoury, Emanuele Fakhoury, and Hendrik Speleers.
\newblock Exsplinet: An interpretable and expressive spline-based neural
  network.
\newblock \emph{Neural Networks}, 152:\penalty0 332--346, 2022.

\bibitem[Feng and Prohl(2003)]{feng2003numerical}
Xiaobing Feng and Andreas Prohl.
\newblock Numerical analysis of the allen-cahn equation and approximation for
  mean curvature flows.
\newblock \emph{Numerische Mathematik}, 94:\penalty0 33--65, 2003.

\bibitem[Ge et~al.(1999)Ge, Hang, and Zhang]{ge1999adaptive}
Shuzhi~Sam Ge, Chang~Chieh Hang, and Tao Zhang.
\newblock Adaptive neural network control of nonlinear systems by state and
  output feedback.
\newblock \emph{IEEE Transactions on Systems, Man, and Cybernetics, Part B
  (Cybernetics)}, 29\penalty0 (6):\penalty0 818--828, 1999.

\bibitem[Girosi and Poggio(1989)]{girosi1989representation}
Federico Girosi and Tomaso Poggio.
\newblock Representation properties of networks: Kolmogorov's theorem is
  irrelevant.
\newblock \emph{Neural Computation}, 1\penalty0 (4):\penalty0 465--469, 1989.

\bibitem[Goswami et~al.(2022)Goswami, Yin, Yu, and
  Karniadakis]{goswami2022physics}
Somdatta Goswami, Minglang Yin, Yue Yu, and George~Em Karniadakis.
\newblock A physics-informed variational deeponet for predicting crack path in
  quasi-brittle materials.
\newblock \emph{Computer Methods in Applied Mechanics and Engineering},
  391:\penalty0 114587, 2022.

\bibitem[Guo and Jiang(2021)]{guo2021construct}
Ruchi Guo and Jiahua Jiang.
\newblock Construct deep neural networks based on direct sampling methods for
  solving electrical impedance tomography.
\newblock \emph{SIAM Journal on Scientific Computing}, 43\penalty0
  (3):\penalty0 B678--B711, 2021.

\bibitem[Harlim et~al.(2021)Harlim, Jiang, Liang, and Yang]{harlim2021machine}
John Harlim, Shixiao~W Jiang, Senwei Liang, and Haizhao Yang.
\newblock Machine learning for prediction with missing dynamics.
\newblock \emph{Journal of Computational Physics}, 428:\penalty0 109922, 2021.

\bibitem[He(2023)]{he2023optimal}
Juncai He.
\newblock On the optimal expressive power of relu dnns and its application in
  approximation with kolmogorov superposition theorem.
\newblock \emph{arXiv preprint arXiv:2308.05509}, 2023.

\bibitem[He et~al.(2015)He, Chen, and Yin]{he2015adaptive}
Wei He, Yuhao Chen, and Zhao Yin.
\newblock Adaptive neural network control of an uncertain robot with full-state
  constraints.
\newblock \emph{IEEE transactions on cybernetics}, 46\penalty0 (3):\penalty0
  620--629, 2015.

\bibitem[Karniadakis et~al.(2021)Karniadakis, Kevrekidis, Lu, Perdikaris, Wang,
  and Yang]{karniadakis2021physics}
George~Em Karniadakis, Ioannis~G Kevrekidis, Lu~Lu, Paris Perdikaris, Sifan
  Wang, and Liu Yang.
\newblock Physics-informed machine learning.
\newblock \emph{Nature Reviews Physics}, 3\penalty0 (6):\penalty0 422--440,
  2021.

\bibitem[Kolmogorov(1957)]{kolmogorov1957representation}
Andrei~Nikolaevich Kolmogorov.
\newblock On the representation of continuous functions of many variables by
  superposition of continuous functions of one variable and addition.
\newblock In \emph{Doklady Akademii Nauk}, volume 114, pages 953--956. Russian
  Academy of Sciences, 1957.

\bibitem[Kolmogorov(1961)]{kolmogorov1961representation}
Andre{{\u\i}}~Nikolaevich Kolmogorov.
\newblock \emph{On the representation of continuous functions of several
  variables by superpositions of continuous functions of a smaller number of
  variables}.
\newblock American Mathematical Society, 1961.

\bibitem[K{\"o}ppen(2002)]{koppen2002training}
Mario K{\"o}ppen.
\newblock On the training of a kolmogorov network.
\newblock In \emph{Artificial Neural Networks ICANN 2002: International
  Conference Madrid, Spain, August 28--30, 2002 Proceedings 12}, pages
  474--479. Springer, 2002.

\bibitem[Kovachki et~al.(2023)Kovachki, Li, Liu, Azizzadenesheli, Bhattacharya,
  Stuart, and Anandkumar]{kovachki2023neural}
Nikola Kovachki, Zongyi Li, Burigede Liu, Kamyar Azizzadenesheli, Kaushik
  Bhattacharya, Andrew Stuart, and Anima Anandkumar.
\newblock Neural operator: Learning maps between function spaces with
  applications to pdes.
\newblock \emph{Journal of Machine Learning Research}, 24\penalty0
  (89):\penalty0 1--97, 2023.

\bibitem[Lai et~al.(2018)Lai, Chang, Yang, and Liu]{lai2018modeling}
Guokun Lai, Wei-Cheng Chang, Yiming Yang, and Hanxiao Liu.
\newblock Modeling long-and short-term temporal patterns with deep neural
  networks.
\newblock In \emph{The 41st international ACM SIGIR conference on research \&
  development in information retrieval}, pages 95--104, 2018.

\bibitem[Lai and Shen(2021)]{lai2021kolmogorov}
Ming-Jun Lai and Zhaiming Shen.
\newblock The kolmogorov superposition theorem can break the curse of
  dimensionality when approximating high dimensional functions.
\newblock \emph{arXiv preprint arXiv:2112.09963}, 2021.

\bibitem[Layton(2008)]{layton2008introduction}
W.~J. Layton.
\newblock \emph{Introduction to the numerical analysis of incompressible
  viscous flows}, volume~6.
\newblock Society for Industrial and Applied Mathematics (SIAM), 2008.

\bibitem[Leni et~al.(2013)Leni, Fougerolle, and Truchetet]{leni2013kolmogorov}
Pierre-Emmanuel Leni, Yohan~D Fougerolle, and Fr{\'e}d{\'e}ric Truchetet.
\newblock The kolmogorov spline network for image processing.
\newblock In \emph{Image Processing: Concepts, Methodologies, Tools, and
  Applications}, pages 54--78. IGI Global, 2013.

\bibitem[Li et~al.(2020)Li, Kovachki, Azizzadenesheli, Liu, Bhattacharya,
  Stuart, and Anandkumar]{li2020fourier}
Zongyi Li, Nikola Kovachki, Kamyar Azizzadenesheli, Burigede Liu, Kaushik
  Bhattacharya, Andrew Stuart, and Anima Anandkumar.
\newblock Fourier neural operator for parametric partial differential
  equations.
\newblock \emph{arXiv preprint arXiv:2010.08895}, 2020.

\bibitem[Liang et~al.(2022)Liang, Huang, and Zhang]{liang2022stiffness}
Senwei Liang, Zhongzhan Huang, and Hong Zhang.
\newblock Stiffness-aware neural network for learning hamiltonian systems.
\newblock In \emph{International Conference on Learning Representations}, 2022.

\bibitem[Lin and Unbehauen(1993)]{lin1993realization}
Ji-Nan Lin and Rolf Unbehauen.
\newblock On the realization of a kolmogorov network.
\newblock \emph{Neural Computation}, 5\penalty0 (1):\penalty0 18--20, 1993.

\bibitem[Liu et~al.(2024{\natexlab{a}})Liu, Ma, Wang, Matusik, and
  Tegmark]{liu2024kan2}
Ziming Liu, Pingchuan Ma, Yixuan Wang, Wojciech Matusik, and Max Tegmark.
\newblock Kan 2.0: Kolmogorov-arnold networks meet science.
\newblock \emph{arXiv preprint arXiv:2408.10205}, 2024{\natexlab{a}}.

\bibitem[Liu et~al.(2024{\natexlab{b}})Liu, Wang, Vaidya, Ruehle, Halverson,
  Solja{\v{c}}i{\'c}, Hou, and Tegmark]{liu2024kan}
Ziming Liu, Yixuan Wang, Sachin Vaidya, Fabian Ruehle, James Halverson, Marin
  Solja{\v{c}}i{\'c}, Thomas~Y Hou, and Max Tegmark.
\newblock Kan: Kolmogorov-arnold networks.
\newblock \emph{arXiv preprint arXiv:2404.19756}, 2024{\natexlab{b}}.

\bibitem[Long et~al.(2019)Long, Lu, and Dong]{long2019pde}
Zichao Long, Yiping Lu, and Bin Dong.
\newblock Pde-net 2.0: Learning pdes from data with a numeric-symbolic hybrid
  deep network.
\newblock \emph{Journal of Computational Physics}, 399:\penalty0 108925, 2019.

\bibitem[Lu et~al.(2021)Lu, Jin, Pang, Zhang, and Karniadakis]{lu2021learning}
Lu~Lu, Pengzhan Jin, Guofei Pang, Zhongqiang Zhang, and George~Em Karniadakis.
\newblock Learning nonlinear operators via deeponet based on the universal
  approximation theorem of operators.
\newblock \emph{Nature machine intelligence}, 3\penalty0 (3):\penalty0
  218--229, 2021.

\bibitem[Meng et~al.(2020)Meng, Li, Zhang, and Karniadakis]{meng2020ppinn}
Xuhui Meng, Zhen Li, Dongkun Zhang, and George~Em Karniadakis.
\newblock Ppinn: Parareal physics-informed neural network for time-dependent
  pdes.
\newblock \emph{Computer Methods in Applied Mechanics and Engineering},
  370:\penalty0 113250, 2020.

\bibitem[Montanelli and Yang(2020)]{montanelli2020error}
Hadrien Montanelli and Haizhao Yang.
\newblock Error bounds for deep relu networks using the kolmogorov--arnold
  superposition theorem.
\newblock \emph{Neural Networks}, 129:\penalty0 1--6, 2020.

\bibitem[Mou et~al.(2021)Mou, Koc, San, Rebholz, and Iliescu]{mou2021data}
Changhong Mou, Birgul Koc, Omer San, Leo~G Rebholz, and Traian Iliescu.
\newblock Data-driven variational multiscale reduced order models.
\newblock \emph{Computer Methods in Applied Mechanics and Engineering},
  373:\penalty0 113470, 2021.

\bibitem[Mou et~al.(2023{\natexlab{a}})Mou, Chen, and
  Iliescu]{mou2023efficient}
Changhong Mou, Nan Chen, and Traian Iliescu.
\newblock An efficient data-driven multiscale stochastic reduced order modeling
  framework for complex systems.
\newblock \emph{Journal of Computational Physics}, 493:\penalty0 112450,
  2023{\natexlab{a}}.

\bibitem[Mou et~al.(2023{\natexlab{b}})Mou, Smith, and Chen]{mou2023combining}
Changhong Mou, Leslie~M Smith, and Nan Chen.
\newblock Combining stochastic parameterized reduced-order models with machine
  learning for data assimilation and uncertainty quantification with partial
  observations.
\newblock \emph{Journal of Advances in Modeling Earth Systems}, 15\penalty0
  (10):\penalty0 e2022MS003597, 2023{\natexlab{b}}.

\bibitem[Pijanowski et~al.(2014)Pijanowski, Tayyebi, Doucette, Pekin, Braun,
  and Plourde]{pijanowski2014big}
Bryan~C Pijanowski, Amin Tayyebi, Jarrod Doucette, Burak~K Pekin, David Braun,
  and James Plourde.
\newblock A big data urban growth simulation at a national scale: configuring
  the gis and neural network based land transformation model to run in a high
  performance computing (hpc) environment.
\newblock \emph{Environmental Modelling \& Software}, 51:\penalty0 250--268,
  2014.

\bibitem[Poggio et~al.(2020)Poggio, Banburski, and Liao]{poggio2020theoretical}
Tomaso Poggio, Andrzej Banburski, and Qianli Liao.
\newblock Theoretical issues in deep networks.
\newblock \emph{Proceedings of the National Academy of Sciences}, 117\penalty0
  (48):\penalty0 30039--30045, 2020.

\bibitem[Popov et~al.(2021)Popov, Mou, Sandu, and
  Iliescu]{popov2021multifidelity}
Andrey~A Popov, Changhong Mou, Adrian Sandu, and Traian Iliescu.
\newblock A multifidelity ensemble kalman filter with reduced order control
  variates.
\newblock \emph{SIAM Journal on Scientific Computing}, 43\penalty0
  (2):\penalty0 A1134--A1162, 2021.

\bibitem[Raissi et~al.(2019)Raissi, Perdikaris, and
  Karniadakis]{raissi2019physics}
Maziar Raissi, Paris Perdikaris, and George~E Karniadakis.
\newblock Physics-informed neural networks: A deep learning framework for
  solving forward and inverse problems involving nonlinear partial differential
  equations.
\newblock \emph{Journal of Computational physics}, 378:\penalty0 686--707,
  2019.

\bibitem[Shen et~al.(1999)Shen, Medjo, and Wang]{shen1999wind}
J.~Shen, T.T. Medjo, and S.~Wang.
\newblock {On a Wind-Driven, Double-Gyre, Quasi-Geostrophic Ocean Model:
  Numerical Simulations and Structural Analysis}.
\newblock \emph{Journal of Computational Physics}, 155\penalty0 (2):\penalty0
  387--409, 1999.
\newblock ISSN 0021-9991.

\bibitem[Shen and Yang(2010)]{shen2010numerical}
Jie Shen and Xiaofeng Yang.
\newblock Numerical approximations of allen-cahn and cahn-hilliard equations.
\newblock \emph{Discrete Contin. Dyn. Syst}, 28\penalty0 (4):\penalty0
  1669--1691, 2010.

\bibitem[Shen et~al.(2018)Shen, Xu, and Yang]{shen2018scalar}
Jie Shen, Jie Xu, and Jiang Yang.
\newblock The scalar auxiliary variable (sav) approach for gradient flows.
\newblock \emph{Journal of Computational Physics}, 353:\penalty0 407--416,
  2018.

\bibitem[Shen et~al.(2019)Shen, Xu, and Yang]{shen2019new}
Jie Shen, Jie Xu, and Jiang Yang.
\newblock A new class of efficient and robust energy stable schemes for
  gradient flows.
\newblock \emph{SIAM Review}, 61\penalty0 (3):\penalty0 474--506, 2019.

\bibitem[Sprecher and Draghici(2002)]{sprecher2002space}
David~A Sprecher and Sorin Draghici.
\newblock Space-filling curves and kolmogorov superposition-based neural
  networks.
\newblock \emph{Neural Networks}, 15\penalty0 (1):\penalty0 57--67, 2002.

\bibitem[Sun and Beckermann(2007)]{sun2007sharp}
Ying Sun and Christoph Beckermann.
\newblock Sharp interface tracking using the phase-field equation.
\newblock \emph{Journal of Computational Physics}, 220\penalty0 (2):\penalty0
  626--653, 2007.

\bibitem[Temam(2024)]{temam2024navier}
Roger Temam.
\newblock \emph{Navier--Stokes equations: theory and numerical analysis},
  volume 343.
\newblock American Mathematical Society, 2024.

\bibitem[Toscano et~al.(2024)Toscano, Oommen, Varghese, Zou, Daryakenari, Wu,
  and Karniadakis]{toscano2024pinns}
Juan~Diego Toscano, Vivek Oommen, Alan~John Varghese, Zongren Zou,
  Nazanin~Ahmadi Daryakenari, Chenxi Wu, and George~Em Karniadakis.
\newblock From pinns to pikans: Recent advances in physics-informed machine
  learning.
\newblock \emph{arXiv preprint arXiv:2410.13228}, 2024.

\bibitem[Trentin(2001)]{trentin2001networks}
Edmondo Trentin.
\newblock Networks with trainable amplitude of activation functions.
\newblock \emph{Neural Networks}, 14\penalty0 (4-5):\penalty0 471--493, 2001.

\bibitem[Winovich et~al.(2019)Winovich, Ramani, and Lin]{winovich2019convpde}
Nick Winovich, Karthik Ramani, and Guang Lin.
\newblock Convpde-uq: Convolutional neural networks with quantified uncertainty
  for heterogeneous elliptic partial differential equations on varied domains.
\newblock \emph{Journal of Computational Physics}, 394:\penalty0 263--279,
  2019.

\bibitem[Wong et~al.(2010)Wong, Xia, and Chu]{wong2010adaptive}
Wai~Keung Wong, Min Xia, and WC~Chu.
\newblock Adaptive neural network model for time-series forecasting.
\newblock \emph{European Journal of Operational Research}, 207\penalty0
  (2):\penalty0 807--816, 2010.

\bibitem[Yan et~al.(2021)Yan, Hao, Zhang, Illman, Lin, and
  Zeng]{yan2021accelerating}
Hengnian Yan, Chenyu Hao, Jiangjiang Zhang, Walter~A Illman, Guang Lin, and
  Lingzao Zeng.
\newblock Accelerating groundwater data assimilation with a gradient-free
  active subspace method.
\newblock \emph{Water Resources Research}, 57\penalty0 (12):\penalty0
  e2021WR029610, 2021.

\bibitem[Zhang et~al.(2024)Zhang, Zhang, Shen, and Lin]{zhang2024energy}
Jiahao Zhang, Shiheng Zhang, Jie Shen, and Guang Lin.
\newblock Energy-dissipative evolutionary deep operator neural networks.
\newblock \emph{Journal of Computational Physics}, 498:\penalty0 112638, 2024.

\bibitem[Zhang et~al.(2022)Zhang, Shen, and Yang]{zhang2022neural}
Shijun Zhang, Zuowei Shen, and Haizhao Yang.
\newblock Neural network architecture beyond width and depth.
\newblock \emph{Advances in Neural Information Processing Systems},
  35:\penalty0 5669--5681, 2022.

\bibitem[Zhu et~al.(2023)Zhu, Feng, Lin, and Lu]{zhu2023fourier}
Min Zhu, Shihang Feng, Youzuo Lin, and Lu~Lu.
\newblock Fourier-deeponet: Fourier-enhanced deep operator networks for full
  waveform inversion with improved accuracy, generalizability, and robustness.
\newblock \emph{Computer Methods in Applied Mechanics and Engineering},
  416:\penalty0 116300, 2023.

\bibitem[Zhuang et~al.(2024)Zhuang, Yao, Zhang, and Karniadakis]{zhuang2024two}
Qiao Zhuang, Chris~Ziyi Yao, Zhongqiang Zhang, and George~Em Karniadakis.
\newblock Two-scale neural networks for partial differential equations with
  small parameters.
\newblock \emph{arXiv preprint arXiv:2402.17232}, 2024.

\end{thebibliography}
